\documentclass[final]{siamltex}
\usepackage{amsfonts}
\usepackage{latexsym,amssymb,amsfonts}
\usepackage{mathrsfs}
\usepackage{graphicx}
\usepackage{psfrag}
\usepackage{amsmath, amsbsy}
\usepackage{amsopn, amstext}
\usepackage{amsmath,multicol,amsopn}
\usepackage{latexsym,amssymb}
\usepackage{amsbsy, amstext,makeidx}
\def\3n{\negthinspace \negthinspace \negthinspace }
\def\2n{\negthinspace \negthinspace }
\def\1n{\negthinspace }

\def\ds{\displaystyle}

\def\={\buildrel \triangle \over =}

\def\resp{{\it resp. }}
\def\mE{{\mathbb{E}}}
\def\b{\beta}

\def\d{\delta}
\def\e{\varepsilon}

\def\l{\lambda}

 \def\n{\nabla}
\def\si{\sigma}
\def\t{\times}
\def\f{\varphi}
\def\th{\theta}
\def\o{\omega}

\def\ns{\noalign{\ss} }

\def\ov{{\overline v}}

\def\pa{\partial}
\def\G{\Gamma}
\def\D{\Delta}

\def\Si{\Sigma}

\def\O{\Omega}

\def\cF{{\cal F}}

\def\cM{{\cal M}}

\def\cP{{\cal P}}

\def\no{\noindent}

\def\ms{\medskip}

\def\q{\quad}
\def\qq{\qquad}

\def\max{\mathop{\rm max}}
\def\min{\mathop{\rm min}}

\def\pa{\partial}

\def\wt{\widetilde}
\def\cd{\cdot}
\def\cds{\cdots}

\def\div{\hbox{\rm div$\,$}}

\def\dist{\hbox{\rm dist$\,$}}

\def\|{\Big |}
\def\({\Big (}
\def\){\Big )}
\def\[{\Big[}
\def\]{\Big]}
\def\be{\begin{equation}}
\def\bel{\begin{equation}\label}
\def\ee{\end{equation}}
\def\bt{\begin{theorem}}
\def\bcd{\begin{condition}}
\def\ecd{\end{condition}}
\def\et{\end{theorem}}
\def\bc{\begin{corollary}}
\def\ec{\end{corollary}}
\def\bde{\begin{definition}}
\def\ede{\end{definition}}
\def\bl{\begin{lemma}}
\def\el{\end{lemma}}
\def\bp{\begin{proposition}}
\def\ep{\end{proposition}}
\def\br{\begin{remark}}
\def\er{\end{remark}}
\def\ba{\begin{array}}
\def\ea{\end{array}}
\def\ed{\end{document}}
\def\ns{\noalign{\ms}}
\def\ds{\displaystyle}

\newtheorem{remark}{Remark}[section]

\newtheorem{condition}{Condition}[section]
\title{OBSERVABILITY ESTIMATE FOR STOCHASTIC SCHR\"{O}DINGER EQUATIONS
AND ITS APPLICATIONS\thanks{This work is
partially supported by the NSF of China under
grant 11101070, and the ERC Advanced Grant
FP7-246775 NUMERIWAVES,  the Grant PI2010-04 of
the Basque Government, the ESF Research
Networking Programme OPTPDE and Grant
MTM2008-03541 of the MICINN, Spain. }}

\author{Qi L\"{u}\thanks{School of Mathematical Sciences, University of Electronic
 Science and Technology of China, Chengdu, 610054, China; and Basque
Center for Applied Mathematics (BCAM),
Mazarredo, 14. 48009 Bilbao Basque Country -
Spain, ({\tt luqi59@163.com}).}}

\begin{document}
\maketitle

\begin{abstract}
In this paper, we establish a boundary
observability estimate for stochastic
Schr\"{o}dinger equations by means of the
global Carleman estimate. Our Carleman estimate
is based on a new fundamental identity for a
stochastic Schr\"{o}dinger-like operator.
Applications to the state observation problem
for semilinear stochastic Schr\"{o}dinger
equations and the unique continuation problem
for stochastic Schr\"{o}dinger equations are
also addressed.
\end{abstract}

\begin{keywords}
stochastic Schr\"{o}dinger equation, global
Carleman estimate, observability estimate,
state observation problem, unique continuation
property
\end{keywords}

\begin{AMS}
93B07, 35B45
\end{AMS}

\pagestyle{myheadings} \thispagestyle{plain}
\markboth{OBSERVABILITY ESTIMATE FOR SSE}{QI
L\"{U}}

\section{Introduction and Main Results}

\q Let $T > 0$, $G \subset \mathbb{R}^{n}$ ($n
\in \mathbb{N}$) be a given bounded domain with
a $C^{2}$ boundary $\G$. Let $\G_0$ be a
suitable chosen nonempty subset (to be given
later) of $\G$. Put $Q \= (0,T) \t G$, $\Si \=
(0,T) \t \G$, and $\Si_0 \= (0,T) \t \G_0$.

Let $(\O, {\cal F}, \{{\cal F}_t\}_{t \geq 0},
P)$ be a complete filtered probability space on
which a  one dimensional standard Brownian
motion $\{ B(t) \}_{t\geq 0}$ is defined. Let
$H$ be a Banach space. Denote by $L^{2}_{\cal
F}(0,T;H)$ the Banach space consisting of all
$H$-valued $\{ {\cal F}_t \}_{t\geq 0}$-adapted
processes $X(\cdot)$ such that
$\mathbb{E}(|X(\cdot)|^2_{L^2(0,T;H)}) <
\infty$; by $L^{\infty}_{\cal F}(0,T;H)$ the
Banach space consisting of all $H$-valued $\{
{\cal F}_t \}_{t\geq 0}$-adapted bounded
processes; and  by $L^{2}_{\cal
F}(\O;C([0,T];H))$ the Banach space consisting
of all $H$-valued $\{ {\cal F}_t \}_{t\geq
0}$-adapted processes $X(\cdot)$ such that
$\mathbb{E}(|X(\cdot)|^2_{C([0,T];H)}) <
\infty$. All of these spaces are endowed with
the canonical norm. Put
$$
H_{T} \= L_{\cal F}^2 (\O; C([0,T];H_{0}^1(G))).
$$

Let us consider the following stochastic
Schr\"{o}dinger equation:
\begin{eqnarray}\label{system1}
\left\{
\begin{array}{lll}
\ds idy + \D ydt = (a_1 \cdot \nabla y + a_2 y + f)dt + (a_3 y + g)dB(t)  &\mbox { in } Q, \\
\ns\ds y = 0 &\mbox{ on } \Si, \\
y(0) = y_0 &\mbox{ in } G,
\end{array}
\right.
\end{eqnarray}
with initial datum $y_0 \in
L^2(\O,\cF_0,P;H_0^1(G))$, suitable
coefficients $a_i$ ($i=1,2,3$), and source
terms $f$ and $g$. The solution to
\eqref{system1} is understood in the following
sense.

\begin{definition}\label{def1}
We call $y\in H_T$ a
solution to the equation \eqref{system1} if \\
1.   $y(0) = y_0$ in $G$, P-a.s.;\\
2. For any $t \in [0,T]$ and $\eta \in
H_0^1(G)$, it holds that
\begin{eqnarray}\nonumber
&\,& \q\int_{G} iy(t,x)\eta(x)dx - \int_{G}
iy(0,x)\eta(x)dx\nonumber\\
&\,& = \int_0^t \int_G \[ \nabla
y(s,x)\cdot\nabla\eta(x) + \big(a_1
\cdot \nabla y + a_2 y + f\big)\eta(x) \]dxds \nonumber \\
&\,&\q + \int_0^t \int_G (a_3 y + g)\eta(x)
dxdB(s), \,\, \mbox { P-a.s. } \nonumber
\end{eqnarray}
\end{definition}

We refer to \cite[Chapter 6]{Prato} for the
well-posedness of the equation \eqref{system1}
in $H_T$, under suitable assumptions (the
assumptions in this paper are enough).

Similar to its deterministic counterpart, the
stochastic Schr\"{o}dinger equation plays an
important role in quantum mechanics. We refer
the readers to \cite{Bar,Kol} and the rich
references therein for the details of its
physical background.

The main purpose of this paper is to establish
a boundary observability estimate for the
equation \eqref{system1} in the following
setting.

Denote by $\nu(x)$ the unit outward normal
vector of $G$ at $x\in \G$.  Let
$x_0\in\big(\mathbb{R}^n\setminus \overline
G\big)$. In what follows, we choose
\begin{equation}\label{G0}
\G_0=\big\{ x\in \G :\, (x-x_0)\cdot \nu(x)>0
\big\}.
\end{equation}
We assume that
\begin{eqnarray}\label{coai}
\left\{\begin{array} {ll}  \ds  i a_1 \in L_{
\mathcal{F}}^{\infty}(0,T;W_0^{1,\infty}(G;\mathbb{R}^{n})),
\\
\ns\ds a_2
\in L_{ \mathcal{F}}^{\infty}(0,T;W^{1,\infty}(G)), \\
\ns  \ds  a_3 \in L_{\mathcal{
F}}^{\infty}(0,T;W^{1,\infty}(G)),
\end{array}
\right.
\end{eqnarray}
and that
\begin{eqnarray}\label{fg}
\left\{
\begin{array}{ll}\ds
f \in L^2_{\mathcal{F}}(0,T;H_0^1(G)), \\
\ns\ds  g \in L^2_{\mathcal{F}}(0,T;H^1(G)).
\end{array}
\right.
\end{eqnarray}

In the sequel, we put
\begin{equation}\label{cA}
r_1\=|a_1|^2_{L_{
\mathcal{F}}^{\infty}(0,T;W_0^{1,\infty}(G;\mathbb{R}^{n}))}
+ |a_2|^2_{L_{
\mathcal{F}}^{\infty}(0,T;W^{1,\infty}(G))} +
|a_3|^2_{L_{
\mathcal{F}}^{\infty}(0,T;W^{1,\infty}(G))} +
1,
\end{equation}
and denote by $C$ a generic positive constant
 depending only on $T$, $G$ and $x_0$, which may change
from line to line.

Now we state the main result of this paper as
follows.

\vspace{0.3cm}

\begin{theorem}\label{observability}
 If the conditions \eqref{G0}--\eqref{fg} hold, then any solution of the equation \eqref{system1} satisfies that
\begin{equation} \label{obser esti2}
\begin{array}{ll}\ds
 \q |y_0|_{L^2(\Omega,{ \mathcal{F}}_0, P; H_0^1(G))} \\
\ns\ds \leq   e^{C r_1}\Big(\Big|\frac{\partial
y}{\partial \nu}\Big |_{L^2_{ \mathcal{
F}}(0,T;L^2(\Gamma_0))} + |f|_{L^2_{ \mathcal{
F}}(0,T;H_0^1(G))} + |g|_{L^2_{ \mathcal{
F}}(0,T;H^1(G))}\Big).
\end{array}
\end{equation}
\end{theorem}

\begin{remark}
Since $y$ belongs only to $H_T$, its normal
derivative $\frac{\pa y}{\pa\nu}$ may not make
sense. Fortunately, due to the hidden
regularity of the solution to the equation
\eqref{system1}, one can show that $\frac{\pa
y}{\pa\nu}$ exists and belongs to
$L^2_{\cF}(0,T;L^2(\G))$(see Proposition
\ref{hregularity} for more details).
\end{remark}

It is well-known that observability estimates
(in the spirit of \eqref{obser esti2}) for
partial differential equations play fundamental
role in proving the controllability of the dual
control systems. There exist many approaches
and results addressing the observability
estimate for determinisitc Schr\"{o}dinger
equations. For example, similar results in the
spirit of Theorem \ref{observability} are
obtained by Carleman estimate (e.g.
\cite{Baudouin-Puel,Lasiecka-Triggiani-Zhang,Mercado-Osses-Rosier}),
by the classical Rellich-type multiplier
approach (\cite{Machtyngier}), by the
microlocal analysis approach
(\cite{Lebeau,Phung}), and so on. We refer to
\cite{Zuazua} for a nice survey in this
respect. However, people know very little about
the stochastic counterpart. To our best
knowledge, \cite{Luqi4} is the only published
result for this problem, where partial results
in this paper have been announced without
detailed proofs.

Besides its important application to the
controllability problem, the observability
estimate not only have its own interest (a kind
of energy estimate and quantitative uniqueness
for the solution) but also has some other
important applications. For instance, a typical
application of this sort of estimates is to
study the state observation problem, that is,
to determine the state of a system by a
suitable observation. Once the observability is
obtained, we may conclude that the state can be
uniquely determined from the observed data and
continuously depends on it. For instance, once
the inequality \eqref{obser esti2} is
established, it follows that $y\in H_T$ is
determined by $\ds\frac{\pa y}{\pa
\nu}\Big|_{(0,T)\times \G_0}$ continuously. In
Section \ref{Sec app}, we shall consider a
state observation problem for semilinear
stochastic Schr\"{o}dinger equations.

In this paper, we will prove Theorem
\ref{observability}  by applying the global
Carleman estimate (See Theorem \ref{thcarleman
est} below).

We now introduce the weight functions to be
used in our Carleman estimate. Let
\begin{equation}\label{psi}
\psi(x) = |x-x_0|^2 + \tau,
\end{equation}
where $\tau$ is a positive constant such that
$\psi \geq \frac{5}{6}|\psi|_{L^{\infty}(G)}$.
Let $s>0$ and $\l>0$.  Put
\begin{equation}\label{lvarphi}
\ell = s\frac{e^{4\l \psi} - e^{5\l
|\psi|_{L^{\infty}(G)}}}{t^2(T-t)^2}, \qq
\varphi = \frac{e^{4\l \psi} }{t^2(T-t)^2},\qq
\theta=e^\ell.
\end{equation}

We have the following global Carleman
inequality.

\vspace{0.2cm}

\begin{theorem}\label{thcarleman est}
According to \eqref{G0}--\eqref{cA} and
\eqref{lvarphi}, there is an $s_1>0$ (depending
on $r_1$) and a $\l_1>0$ such that for each
$s\geq s_1$, $\l\geq \l_1$ and for any solution
of the equation \eqref{system1}, it holds that
\begin{eqnarray}\label{carleman est}
\begin{array}{ll}
\ds \q\mathbb{E}\int_Q
\theta^2\Big(s^3\l^4\varphi^3 |y|^2 +
s\l\varphi
|\nabla y|^2\Big) dxdt \\
\ns \ds \leq  C \Big\{\mathbb{E}\int_Q \theta^2
\Big(|f|^2 +
 s^2\l^2\varphi^2 |g|^2 + |\nabla g|^2 \Big)dxdt +
\mathbb{E}\int_0^T\int_{\G_0}\theta^2
s\l\varphi\Big| \frac{\pa y}{\pa \nu}\Big|^2d\G
dt \Big\}.
\end{array}
\end{eqnarray}
Further, if $g\in
L^2_\cF(0,T;H^1(G;\mathbb{R}))$, then
\eqref{carleman est} can be strengthened as the
following:
\begin{eqnarray}\label{carleman est1}
\begin{array}
{ll} \ds \q\mathbb{E}\int_Q
\theta^2\Big(s^3\l^4\varphi^3 |y|^2 +
s\l\varphi
|\nabla y|^2\Big) dxdt \\
\ns \ds\leq    C \Big\{\mathbb{E}\int_Q
\theta^2 \Big(|f|^2 +
 s^2\l^2\varphi^2 |g|^2  \Big)dxdt +
\mathbb{E}\int_0^T\int_{\G_0}\theta^2
s\l\varphi\Big| \frac{\pa y}{\pa \nu}\Big|^2d\G
dt \Big\}.
\end{array}
\end{eqnarray}

\end{theorem}

Carleman estimate is an important tool for the
study of unique continuation property,
stabilization, controllability and inverse
problems for deterministic partial differential
equations (e.g.
\cite{Baudouin-Puel,Lasiecka-Triggiani-Zhang,Mercado-Osses-Rosier,Yamamoto,Zhangxu1,Zuazua}).
Although there are numerous results for the
Carleman estimate for deterministic partial
differential equations, people know very little
about the corresponding stochastic situation.
In fact, as far as we know,
\cite{barbu1,Luqi4,Luqi5,Tang-Zhang1,Zhangxu3}
are the only five published papers addressing
the Carleman estimate for stochastic partial
differential equations. The references
\cite{barbu1,Luqi5,Tang-Zhang1} are devoted to
stochastic heat equations, while
\cite{Zhangxu3} is concerned with stochastic
wave equations. In \cite{Luqi4},   Theorem
\ref{thcarleman est} was announced without
proof.

At first glance, the proof of Theorem
\ref{thcarleman est} looks very similar to that
of the global Carleman estimate for
(stochastic) parabolic equations (See
\cite{Fursikov-Imanuvilov1,Tang-Zhang1}).
Furthermore, one can find that the idea behind
the proofs in this paper and
\cite{Fursikov-Imanuvilov1,Tang-Zhang1} are
analogous. Nevertheless, the specific proofs
have big differences. First, we have to choose
different weight functions. Second, we deal
with different equations. Such kind of
differences lead to considerably different
difficulties in the proof of Theorem
\ref{thcarleman est}. One cannot simply mimic
the proofs in
\cite{Fursikov-Imanuvilov1,Tang-Zhang1} to
obtain Theorem \ref{thcarleman est}. Indeed,
even in the deterministic setting, the proof of
the global Carleman estimate for
Schr\"{o}dinger equations are much more
complicated than that for the parabolic and
hyperbolic equations (see
\cite{Zhangxu5,Lasiecka-Triggiani-Zhang}).

The rest of this paper is organized as follows.
In Section 2, we
 give some preliminary results, including an energy estimate
and the hidden regularity for solutions of the
equation \eqref{system1}. Section 3 is
addressed to establish a crucial identity for a
stochastic Schr\"{o}dinger-like operator. Then,
in Section 4, we derive the desired Carleman
estimate.  Section 5 is devoted to prove
Theorem \ref{observability}. In Section
\ref{Sec app}, as applications of the
observability/Carleman estimates developed in
this work, we study a state observation problem
for semilinear stochastic Schr\"{o}dinger
equations and establish a unique continuation
property for the solution to the equation
\eqref{system1}. Finally, we present some
further comments and open problems concerned
with this paper in Section 7.

\section{Some preliminaries}

In this section, we give some preliminary
results which will be used later.

To begin with, for the sake of completeness, we
give an energy estimate for the equation
\eqref{system1}.

\vspace{0.1cm}

\begin{proposition} \label{Oprop1}
According to \eqref{G0}--\eqref{cA},  for all
$y$ which solve the equation \eqref{system1},
it holds that
\begin{equation}\label{energyesti1}
\mathbb{E}|y(t)|^2_{H_0^1(G)} \leq e^{C r_1}
\Big(\mathbb{E}|y(s)|^2_{H^1_0(G)} +
|f|^2_{L^2_{\mathcal{F}}(0,T;H^1_0(G))} +
|g|^2_{L^2_{\mathcal{F}}(0,T;H^1_0(G))}\Big),
\end{equation}
for any $ s, t\in [0, T]$.
\end{proposition}

\vspace{0.1cm}

{\em Proof }:  Without loss of generality, we
assume that $t <s$. To begin with, we compute $
\mathbb{E}| y(t)|^2_{ L^2(G)} - \mathbb{E}|
y(s)|^2_{ L^2(G)}$  and $ \mathbb{E}|\nabla
y(t)|^2_{ L^2(G)} - \mathbb{E}|\nabla y(s)|^2_{
L^2(G)}$. The first one reads
\begin{equation}\label{Eyt}
\begin{array}{ll} \ds
\mathbb{E}| y(t)|^2_{ L^2(G)} - \mathbb{E}| y(s)|^2_{
L^2(G)}\\
\ns  \ds =  -\mathbb{E}\int_t^s\int_G \big(y
d\bar{y}+\bar{y}dy + dy
d\bar{y}\big)dx\\
\ns\ds  =  \mathbb{E}\int_t^s\int_G \Big\{ i
y\big(\D \bar{y} - a_1\cdot \nabla\bar{y} - a_2
\bar{y} -\bar{f}\big) - i\bar{y}\big(\D y -
a_1\cdot\nabla y - a_2 y - f\big) \\
\ns \ds\q - \big(a_3 y + g\big)\big(a_3 \bar{y} + \bar{g}\big) \Big\}dxd\si \\
\ns\ds   =   \mathbb{E}\int_t^s\int_G \Big\{ i
\big[\div(y\nabla\bar{y})-|\nabla y|^2 - \div(|
y|^2 a_1) +
\div(a_1)|y|^2 - a_2|y|^2 - y\bar{f}\, \big] \\
\ns \ds \q - i \big[\div(\bar{y}\nabla
y)-|\nabla y|^2 - \div(| y|^2 a_1) +
\div(a_1)|y|^2 - a_2|y|^2 - f\bar{y} \big] \\
\ns \ds\q  - (a_3 y + g)(a_3 \bar{y} + \bar{g}) \Big\}dxd\si\\
\ns\ds  \leq \mathbb{E}\int_t^s
2\Big[\big(|a_3|_{L^{\infty}(G)}+1\big)|y|^2_{L^2(G)}+
|f|_{L^2(G)}^2 + |g|^2_{L^2(G)}\Big]dxd\si.
\end{array}
\end{equation}
The second one is
\begin{equation}\label{Etyt}
\begin{array}{ll}
\q\mathbb{E}|\nabla y(t)|^2_{ L^2(G)} -
\mathbb{E}|\nabla
y(s)|^2_{L^2(G)}\\
\ns\ds  = -\mathbb{E}\int_t^s\int_G \big(\nabla
y d\n\bar{y} + \nabla
\bar{y} d\n y + d\nabla y d\nabla \bar{y}\big)dx\\
\ns  \ds =  -\mathbb{E}\int_t^s\int_G \Big\{
\div(\nabla y d\bar{y}) - \D y d\bar{y} +
\div(\nabla\bar{y} dy) - \D \bar{y}dy +
d\nabla y d\nabla \bar{y} \Big\}dx\\
\ns \ds  =  -\mathbb{E}\int_t^s\int_G \Big\{\D
y \Big[i\big(\D\bar{y} - a_1\cdot\nabla \bar{y}
-a_2\bar{y} -f \big)\Big]-\D\bar{y}\Big[
i\big(\D y -
a_1\cdot\nabla y - a_2 y - f\big)\Big]\\
\ns \ds  \q  +\nabla(a_3 y +
g)\nabla(a_3\bar{y}+\bar{g})
\Big\}dxd\si \\
\ns  \ds  \leq  2\mathbb{E}\int_t^s
\Big\{\big(|a_1|^2_{W^{1,\infty}(G;\mathbb{R}^m)}+|a_3|^2_{W^{1,\infty}(G)}+1\big)|\nabla
y|^2_{L^2(G)}\\
\ns\ds  \q  +
\big(|a_2|^2_{W^{1,\infty}(G)}+|a_3|^2_{W^{1,\infty}(G)}+1\big)|y|^2_{L^2(G)}
+|f|^2_{H_0^1(G)} +
|g|^2_{H^1_0(G)}\Big\}dxd\si.
\end{array}
\end{equation}
From \eqref{Eyt} and \eqref{Etyt}, we have that
\begin{equation}\label{energyest2}
\begin{array}{ll}\ds
 \q\mathbb{E}| y(t)|^2_{H_0^1(G)} -
\mathbb{E}|y(s)|^2_{H_0^1(G)} \\
\ns\ds \leq 2(r_1+1)\mathbb{E}\int_t^s\int_G
|y(\si)|^2_{H_0^1(G)}dxd\si +
\mathbb{E}\int_t^s\int_G
\big(|f(\si)|^2_{H_0^1(G)}+|g(\si)|^2_{H^1_0(G)}\big)dxd\si.
\end{array}
\end{equation}
From this, and thanks to Gronwall's inequality,
we arrive at
\begin{equation}\label{energyest3}
\mathbb{E}| y(t)|^2_{H_0^1(G)}\leq
e^{2(r_1+1)}\Big\{\mathbb{E}|y(s)|^2_{H_0^1(G)}
+ \mathbb{E}\int_0^T\int_G
\big(|f|^2_{H_0^1(G)}+|g|^2_{H^1_0(G)}\big)dxd\si\Big\},
\end{equation}
which implies the inequality
\eqref{energyesti1} immediately.
\endproof

\vspace{0.1cm}

\begin{remark}
The proof of this proposition is almost
standard.  However, people may doubt the
correctness of the inequality
\eqref{energyesti1} for $t<s$ because of the
very fact that the equation \eqref{system1} is
time irreversible. Fortunately, the inequality
\eqref{energyesti1} is true for $t<s$. In fact,
in the stochastic setting one should divide the
time irreversible systems into two classes. The
first class of time irreversibility is caused
by the energy dissipation. Thus, one cannot
estimate the energy of the system at time $t$
by that at time $s$ uniformly when $t<s$. A
typical example of such kind of systems is the
heat equation. The second class of time
irreversibility comes from the stochastic
noise. Such kind of system cannot be solved
backward, that is, if we give the final data
rather than the initial data, then the system
is not well-posed (Recall that, this is the
very starting point of backward stochastic
differential equations). Stochastic
Schr\"{o}dinger equations and stochastic wave
equations are typical systems of the second
class. For these systems, we  can still
estimate the energy at time $t$ by that at time
$s$ for $t<s$.

\end{remark}

Next, we give a  result concerning the hidden
regularity for solutions of the equation
\eqref{system1}. It shows that, solutions of
this equation enjoy a higher regularity on the
boundary than the one provided by the classical
trace theorem for Sobolev spaces.
\vspace{0.2cm}
\begin{proposition}\label{hregularity}
According to \eqref{G0}--\,\eqref{cA}, for any
solution of the equation \eqref{system1}, it
holds that
\begin{equation} \label{hregularity1}
\begin{array}{ll}\ds
\q\Big|\frac{\partial y}{\partial \nu}\Big
|^2_{L^2_{ \mathcal{ F}}(0,T;L^2(\Gamma_0))}\\
\ns\ds \leq  e^{C r_1 }
\Big(|y_0|^2_{L^2(\Omega,{ \mathcal{F}}_0, P;
H_0^1(G))} +|f|^2_{L^2_{ \mathcal{
F}}(0,T;H_0^1(G))} +  |g|^2_{L^2_{ \mathcal{
F}}(0,T;H^1(G))}\Big).
\end{array}
\end{equation}
\end{proposition}
\vspace{0.1cm}

\begin{remark}\label{rm2}
By means of Proposition \ref{hregularity}, we
know that $\ds\Big|\frac{\partial y}{\partial
\nu}\Big |^2_{L^2_{ \mathcal{
F}}(0,T;L^2(\Gamma_0))}$ makes sense. Compared
with Theorem \ref{observability}, Proposition
\ref{hregularity} tells us the fact that
$\ds\Big|\frac{\partial y}{\partial \nu}\Big
|^2_{L^2_{ \mathcal{ F}}(0,T;L^2(\Gamma_0))}$
can be bounded by the initial datum and
non-homogenous terms. This result is the
converse of Theorem \ref{observability} in some
sense.
\end{remark}

To prove Proposition \ref{hregularity}, we
first establish a pointwise identity. For
simplicity, here and in the sequel, we adopt
the notation $\ds y_i \equiv y_{i}(x) \=
\frac{\partial y(x)}{\partial x_i}$, where
$x_i$ is the $i$-th coordinate of a generic
point $x=(x_1,\cdots, x_n)$ in
$\mathbb{R}^{n}$. In a similar manner, we  use
the notation $z_i$, $v_i$, etc., for the
partial derivatives of $z$ and $v$ with respect
to $x_i$.

\medskip

\begin{proposition}\label{prop2}
Let $\mu = \mu(x) =
(\mu^1,\cdots,\mu^n):\mathbb{R}^n \to
\mathbb{R}^n$ be a vector field of class $C^1$
and $z$ an $H^2_{loc}(\mathbb{R}^n)$-valued
$\{\mathcal{F}_t\}_{t\geq 0}$-adapted process.
Then for a.e. $x \in \mathbb{R}^n$ and P-a.s.
$\omega \in \Omega$, it holds that
\begin{eqnarray}\label{identity2}
\begin{array}
{ll} & \ds \mu\cdot\nabla\bar{z}(i dz + \Delta
z dt) +
\mu\cdot\nabla z(-i d\bar{z} + \Delta \bar{z} dt)\\
\ns  =& \ds  \nabla\cd \Big[ (\mu\cdot\nabla
\bar{z})\nabla z+ (\mu\cdot\nabla z)\nabla
\bar{z}  - i (z d\bar{z}) \mu - |\nabla
z|^2\mu \Big]dt + d(i\mu\cd\nabla \bar{z} z)\\
\ns &  \ds - 2\sum_{j,k=1}^n \mu^k_j
z_{j}\bar{z}_{k}dt + (\nabla\cdot \mu) |\nabla
z|^2 dt + i(\nabla\cdot \mu) z d\bar{z} -
i(\mu\cd\nabla d\bar z) dz.
\end{array}
\end{eqnarray}
\end{proposition}

\medskip

 {\em Proof of
Proposition \ref{prop2}} : The proof is a
direct computation. We have that
\begin{eqnarray}\label{h1}
\begin{array}
{ll} & \ds \sum_{k=1}^n\sum_{j=1}^n
\mu^k\bar{z}_k z_{jj}+
\sum_{k=1}^n\sum_{j=1}^n \mu^k z_k \bar{z}_{jj}\\
\ns = & \ds
\sum_{k=1}^n\sum_{j=1}^n\Big[(\mu^k\bar{z}_k
z_j)_j + (\mu^k z_k\bar{z}_j)_j +
\mu^k_k|z_j|^2 - (\mu^k|z_j|^2)_k - 2\mu^k_j
\bar{z}_k z_j \Big]
\end{array}
\end{eqnarray}
and that
\begin{equation}\label{h2}
\begin{array}{ll}\ds
\q i\sum_{k=1}^n(\mu^k\bar{z}_k dz-\mu^k z_k d\bar{z})\\
\ns\ds = i\sum_{k=1}^n\Big[\,d(\mu^k\bar{z}_k
z) - \mu^k z d\bar z_k - \mu^k d\bar{z}_k
dz -(\mu^k zd\bar{z})_k + \mu^k z d\bar z_k + \mu_k^k z d\bar{z}\, \Big]\\
\ns\ds = i\sum_{k=1}^n\Big[\,d(\mu^k\bar{z}_k
z) -  \mu^k d\bar{z}_k dz -(\mu^k zd\bar{z})_k
+ \mu_k^k z d\bar{z} \,\Big].
\end{array}
\end{equation}
Combining \eqref{h1} and \eqref{h2}, we get the
equality \eqref{identity2}.
\endproof

\vspace{0.2cm}

By virtue of  Proposition \ref{prop2}, the
proof of Proposition \ref{hregularity} is
standard. We only give a sketch here.

\vspace{0.2cm}

 {\em Sketch of the Proof of Proposition \ref{hregularity}} : Since $\Gamma $ is $ C^2$,
one can find a vector field $\mu_0 = (\mu_0^1,
\cdots, \mu_0^n) \in
C^1(\overline{G};\mathbb{R}^n)$ such that
$\mu_0 = \nu$ on $\Gamma$(see \cite[page
18]{Komornik} for the construction of $\mu_0$).
Letting $\mu = \mu_0$ and $z = y$ in
Proposition \ref{prop2}, integrating it in $Q$
and taking the expectation, by means of
Proposition \ref{prop2}, with similar
computation in \cite{Zhangxu1}, Proposition
\ref{hregularity} can be obtained immediately.

\section{An Identity for a Stochastic Schr\"{o}dinger-like Operator}

\q  In this section, we obtain an identity for
a stochastic schr\"{o}dinger-like operator,
which is similar to the formula
\eqref{identity2} in the spirit but it takes a
more complex form and play a key role in the
proof of Theorem \ref{thcarleman est}.

\vspace{0.2cm}

Let $\b(t,x)\in
C^{2}(\mathbb{R}^{1+n};\mathbb{R})$, and let
$b^{jk}(t,x)\in
C^{1,2}(\mathbb{R}^{1+n};\;\mathbb{R})$ satisfy
\begin{equation}\label{bjk}
b^{jk}=b^{kj},\qq j,k=1,2,\cdots,n.
\end{equation}
Let us define a (formal) second order
stochastic partial differential operator $\cP$
as
\begin{equation}\label{cp}
\ds
 \cP z \= i\b(t,x)dz+\sum_{j,k=1}^n(b^{jk}(t,x)z_j)_k dt,
 \q i=\sqrt{-1}.
\end{equation}
We have the following equality concerning
$\cP$:

\vspace{0.1cm}
\begin{theorem}\label{identity1}
Let $\ell,\;\Psi\in
C^2(\mathbb{R}^{1+n};\;\mathbb{R})$. Assume
that $z$ is an
$H^2_{loc}(\mathbb{R}^n,\mathbb{C})$-valued
$\{\cF_t\}_{t\geq 0}$-adapted process.
 Put $ v=\th
z$(recall \eqref{lvarphi} for the definition of
$\th$). Then for a.e. $x \in \mathbb{R}^n$ and
P-a.s. $\o\in \O$,  it holds that

\begin{eqnarray}\label{c2a1}
\begin{array}{ll}
&\th(\cP z\overline {I_1}+\overline{\cP z} I_1)+dM+\div V \\
\ns  = & \ds 2|I_1|^2 dt +\sum_{j,k=1}^n
 c^{jk}(v_k\ov_j+\ov_k v_j)
dt +D|v|^2 dt \\
\ns & \ds  +i\sum_{j,k=1}^n\[(\b
b^{jk}\ell_j)_t+
b^{jk}(\b\ell_t)_j\](\ov_kv-v_k\ov) dt
\\
\ns & \ds +i\[\b\Psi+\sum_{j,k=1}^n(\b
b^{jk}\ell_j)_k\](\ov
dv-vd\ov)\\
& \ds  + (\b^2\ell_t)dvd\ov + i\sum_{j,k=1}^n
\b b^{jk}\ell_j (dv d\ov_k - dv_kd\ov),
\end{array}
\end{eqnarray}
where
\begin{equation}\label{c2a2}
\left\{
\begin{array}{ll}\ds I_1\= - i\b \ell_t
v - 2\sum_{j,k=1}^n b^{jk}\ell_j v_k +
\Psi v, \\
\ns\ds A\=\sum_{j,k=1}^n
b^{jk}\ell_j\ell_k-\sum_{j,k=1}^n(b^{jk}\ell_j)_k
-\Psi,
\end{array}
\right.
\end{equation}

\begin{equation}
\label{c2a3} \left\{
\begin{array}{ll}\ds
M\=\b^2\ell_t |v|^2 + i\b\sum_{j,k=1}^n b^{jk}\ell_j(\ov_kv-v_k\ov),\\
\ns\ds V\=[V^1,\cdots,V^k,\cdots,V^n],\\
 \ns\ds V^k\=-i \b\sum_{j=1}^n\[b^{jk}\ell_j(vd\ov -\ov
 dv ) + b^{jk}\ell_t(v_j\ov-\ov_jv) dt\]\\
\ns\ds\qq\,\,\,\, - \Psi\sum_{j=1}^n b^{jk}(v_j\ov+\ov_jv) dt + \sum_{j=1}^n b^{jk}(2A\ell_j+\Psi_j)|v|^2 dt \\
\ns\ds\qq\q+\sum_{j,j',k'=1}^n\(2b^{jk'}b^{j'k}-b^{jk}b^{j'k'}\)\ell_j(v_{j'}\ov_{k'}+\ov_{j'}v_{k'})
dt,
\end{array}
\right.
\end{equation}
and
\begin{equation}\label{cc2a3}
\left\{
\begin{array}{ll}\ds
c^{jk}\= \sum_{j',k'=1}^n\[2(b^{j'k}\ell_{j'})_{k'}b^{jk'}-(b^{jk}b^{j'k'}\ell_{j'})_{k'}\] - b^{jk}\Psi,\\
\ns\ds D\=(\b^2\ell_t)_t
+\sum_{j,k=1}^n(b^{jk}\Psi_k)_j+2\[\sum_{j,k=1}^n(b^{jk}\ell_jA)_k+A\Psi\].
\end{array}
\right.
\end{equation}

\end{theorem}
\medskip

\begin{remark}
Since we only assume that $(b^{jk})_{1\leq
j,k\leq n}$ is symmetric and do not assume that
it is positively definite, then similar to
\cite{Fu} and based on the identity
\eqref{c2a1} in Theorem \ref{identity1}, we can
deduce observability estimate not only for the
stochastic Schr\"odinger equation, but also for
deterministic hyperbolic, Schr\"{o}dinger and
plate equations, which had been derived via
Carleman estimate (see \cite{FYZ},
\cite{Lasiecka-Triggiani-Zhang} and
\cite{Zhangxu5}, respectively).
\end{remark}

 {\em Proof of Theorem \ref{identity1}}: The proof is divided
into three steps.

{\bf Step 1.} By the definition of $v$ and $w$,
a straightforward computation shows that:
\begin{eqnarray}\label{th1eq1}
\begin{array}
{ll} \ds \theta \cP z &=  \ds  i\b dv -   i\b
\ell_t v dt + \sum_{j,k=1}^n
(b^{jk}v_j)_k dt \\
\ns & \ds \q + \sum_{j,k=1}^n b^{jk}\ell_j
\ell_k v dt - 2\sum_{j,k=1}^n b^{jk}\ell_j v_k
dt - \sum_{j,k=1}^n (b^{jk}\ell_j)_k v dt
\\
\ns &=  \ds I_1dt + I_2,
\end{array}
\end{eqnarray}
where
\begin{equation}\label{I2}
I_2 = i\b dv + \sum_{j,k=1}^n (b^{jk}v_j)_kdt +
Avdt.
\end{equation}
Hence we obtain that
\begin{equation}\label{th1eq2}
\theta (Pz\overline{I_1} + \overline{Pz}I_1) =
2|I_1|^2dt + (I_1\overline{I_2} +
I_2\overline{I_1}).
\end{equation}

{\bf Step 2.} In this step, we compute
$I_1\overline{I_2} + I_2\overline{I_1}$. Denote
the three terms in $I_1$ and $I_2$ by $I_1^j$
and $I_2^j$, $j = 1,2,3$, respectively. Then we
have that
\begin{equation}\label{th1eq3}
\begin{array}{ll}\ds
\q    I_1^1\overline{I_2^1}  + I_2^1 \overline{I_1^1}  \\
\ns\ds  = -i\b \ell_t v \overline{ (i\b dv)} + i\b dv \overline{ (-i\b \ell_t v)} \\
\ns\ds = -d(\b^2 \ell_t |v|^2) + (\b^2
\ell_t)_t|v|^2dt + \b^2 \ell_t dvd\ov.
\end{array}
\end{equation}
Noting that
\begin{eqnarray}\label{th1eq4}
\left\{
\begin{array}{lll}
\ds 2vd\ov = d(|v|^2) - (\ov dv - vd\ov) - dv d\ov,\\
\ns\ds 2v \ov_k = (|v|^2)_k - (\ov v_k -
v\ov_k),
\end{array}
\right.
\end{eqnarray}
we find first
\begin{equation}\label{th1eq4.1}
\begin{array}{ll}
\ds & \ds 2i\sum_{j,k=1}^n (\b b^{jk}\ell_j vd\ov)_k \\
\ns=& \ds i \sum_{j,k=1}^n \Big\{\b b^{jk}\ell_j \[  d(|v|^2) - (\ov dv - vd\ov) - dv d\ov \] \Big\}_k \\
\ns  = & \ds i  \sum_{j,k=1}^n\Big\{ \big(\b
b^{jk}\ell_j\big)_k d(|v|^2)
+ \b b^{jk}\ell_j\big[d(|v|^2)\big]_k -  \big[\b b^{jk}\ell_j (\ov dv - vd\ov)\big]_k\\
\ns& \ds \qq\q - \big( \b b^{jk}\ell_j\big)_k
dv d\ov - \b b^{jk}\ell_jdv_k d\ov - \b
b^{jk}\ell_jdvd\ov_k  \Big\},
\end{array}
\end{equation}
next
\begin{equation}\label{th1eq4.2}
\begin{array}{ll}
\ds \q -2i\sum_{j,k=1}^n \big(\b b^{jk}\ell_j \big)_k vd\ov \\
\ns\ds= -i \sum_{j,k=1}^n  \big(\b b^{jk}\ell_j \big)_k \[  d(|v|^2) - (\ov dv - vd\ov) - dv d\ov \] \\
\ns\ds  = -i  \sum_{j,k=1}^n\[ \big(\b
b^{jk}\ell_j\big)_k d(|v|^2)   - \big(\b
b^{jk}\ell_j \big)_k(\ov dv - vd\ov) \ - \big(
\b b^{jk}\ell_j \big)_k dv d\bar v  \],
\end{array}
\end{equation}
then
\begin{equation}\label{th1eq4.3}
\begin{array}{ll}\ds
\q -2i \sum_{j,k=1}^n d\big(\b b^{jk}\ell_j v\ov_k \big) \\
\ns\ds = - i \sum_{j,k=1}^n d\Big\{ \b b^{jk}\ell_j \big[(|v|^2)_k - (\ov v_k - v\ov_k)  \big]\Big\}\\
\ns\ds = - i \sum_{j,k=1}^n \Big\{\big(\b
b^{jk}\ell_j \big)_t (|v|^2)_kdt + \b
b^{jk}\ell_j d\big[(|v|^2)_k\big] - d\big[ \b
b^{jk}\ell_j (\ov v_k - v\ov_k) \big] \Big\},
\end{array}
\end{equation}
and that
\begin{equation}\label{th1eq4.4}
\begin{array}{ll}\ds
\q 2i \sum_{j,k=1}^n  \big(\b b^{jk}\ell_j  \big)_t v\ov_k dt \\
\ns\ds =  i \sum_{j,k=1}^n d \big(\b b^{jk}\ell_j  \big)_t \big[(|v|^2)_k - (\ov v_k - v\ov_k)  \big]dt \\
\ns\ds =  i \[\sum_{j,k=1}^n  \big(\b
b^{jk}\ell_j \big)_t (|v|^2)_kdt   -  \big(\b
b^{jk}\ell_j \big)_t (\ov v_k - v\ov_k)dt
\].
\end{array}
\end{equation}
From \eqref{th1eq4.1}--\eqref{th1eq4.4}, we get
that
\begin{equation}\label{th1eq5}
\begin{array}{ll}\ds
\q (I_1^2 + I_1^3)\overline{I_2^1} +
I_2^1(\overline{I_1^2} + \overline{I_1^3})
\\
\ns\ds = \(- 2\sum_{j,k=1}^n b^{jk}\ell_j v_k +
\Psi v\) \overline{ (i\b dv) } +  i\b dv
\overline{ \( - 2\sum_{j,k=1}^n b^{jk}\ell_j
v_k +
\Psi v \) }\\
\ns\ds = 2i\sum_{j,k=1}^n \b b^{jk}\ell_j  (v_k
d\bar v  - \bar v_k dv) + i\b\Psi (\bar v dv -
vd\bar v)
\\
\ns\ds
 = 2i\sum_{j,k=1}^n \[ \big(\b b^{jk}\ell_j vd\ov\big)_k - \big(\b
b^{jk}\ell_j\big)_k vd\ov - \b b^{jk}\ell_j
vd\ov_k
\]
\\ \ns\ds \q -2i \sum_{j,k=1}^n
\[ d\big(\b b^{jk}\ell_j v\ov_k\big) - \big(\b
b^{jk}\ell_j\big)_t v\ov_k dt - \b b^{jk}\ell_j
vd\ov_k \]
\end{array}
\end{equation}
\begin{equation}
\begin{array}{ll}
\ds\q + 2i\sum_{j,k=1}^n \b b^{jk}\ell_j dv
d\ov_k + i\b\Psi (\bar v dv - vd\bar v)\nonumber\\
 \ns\ds = -i \sum_{j,k=1}^n \[ \b b^{jk}\ell_j(\ov dv
- vd\ov )  \]_k dt -i \sum_{j,k=1}^n d\[ \b
b^{jk}\ell_j(v\ov_k - \ov v_k)  \]  \\
\ns\ds \q - i\sum_{j,k=1}^n (\b b^{jk}\ell_j)_t
(\ov v_k - v \ov_k)dt + i\[ \b\Psi +
\sum_{j,k=1}^n
(\b b^{jk}\ell_j)_k \](\ov dv - vd\ov) \\
\ns\ds \q + i\sum_{j,k=1}^n \b b^{jk}\ell_j (dv
d\ov_k - dv_kd\ov).
\end{array}
\end{equation}

Noting that $b^{jk} = b^{kj}$, we have that
\begin{equation}\label{th1eq7}
\begin{array}{ll}\ds
 \q I_1^1\overline{I_2^2} + I_2^2  \overline{I_1^1}
 \\
 \ns\ds = -i\b \ell_t v \overline{\sum_{j,k=1}^n(b^{jk}v_j)_k}dt + \sum_{j,k=1}^n(b^{jk}v_j)_k \overline{(-i\b \ell_t v)} \\
\ns\ds  = \sum_{j,k=1}^n \[ i\b b^{jk}\ell_t
(v_j \ov - \ov_j v)\]_k dt + i \sum_{j,k=1}^n
b^{jk}(\b \ell_t)_k (\ov_j v - v_j\ov)dt.
\end{array}
\end{equation}
Utilizing $b^{jk} = b^{kj}$ once more, we find
$$
\sum_{j,k,j',k'=1}^n b^{jk}b^{j'k'}\ell_j
(v_{j'}\ov_{kk'} +
\ov_{j'}v_{kk'})=\sum_{j,k,j',k'=1}^n
b^{jk}b^{j'k'}\ell_j (v_{j'k}\ov_{k'} +
\ov_{j'k}v_{k'}).
$$
Hence, we obtain that
\begin{equation}\label{th1eq8}
\begin{array}{ll}\ds
 \q 2\sum_{j,k,j',k'=1}^n b^{jk}b^{j'k'}\ell_j (v_{j'}\ov_{kk'} +
\ov_{j'}v_{kk'})dt \\
\ns\ds = \sum_{j,k,j',k'=1}^n
b^{jk}b^{j'k'}\ell_j (v_{j'}\ov_{kk'} +
\ov_{j'}v_{kk'})dt + \sum_{j,k,j',k'=1}^n
b^{jk}b^{j'k'}\ell_j (v_{j'k}\ov_{k'} +
\ov_{j'k}v_{k'})dt\\
 \ns\ds =
\sum_{j,k,j',k'=1}^n b^{jk}b^{j'k'}\ell_j
(v_{j'}\ov_{k'} +
\ov_{j'}v_{k'})_k dt \\
\ns\ds = \!\sum_{j,k,j',k'=1}^n \!\[
b^{jk}b^{j'k'}\ell_j (v_{j'}\ov_{k'} +
\ov_{j'}v_{k'}) \]_kdt -
\!\!\sum_{j,k,j',k'=1}^n
(b^{jk}b^{j'k'}\ell_j)_k (v_{j'}\ov_{k'} + \ov_{j'}v_{k'})dt.\\
\end{array}
\end{equation}

By the equality \eqref{th1eq8}, we get that

\medskip
\begin{equation}\label{th1eq9}
\begin{array}{ll}\ds
\q I_1^2\overline{I_2^2}   +  I_2^2
\overline{I_1^2}\\
\ns\ds = - 2\!\sum_{j,k=1}^n b^{jk}\ell_j v_k\!
\overline{ \sum_{j,k=1}^n(b^{jk}v_j)_k }dt -  2\!\sum_{j,k=1}^n(b^{jk}v_j)_k  \overline{\sum_{j,k=1}^n b^{jk}\ell_j v_k}dt\\
\ns\ds =  - 2\!\! \sum_{j,k,j',k'=1}^n \!\!\[
b^{jk}b^{j'k'}\ell_j (v_{j'}\ov_{k}\! +\!
\ov_{j'}v_{k}) \]_{k'} dt \!+\!
2\!\!\sum_{j,k,j',k'=1}^n\!\!
b^{j'k'}(b^{jk}\ell_j)_{k'} (v_{j'}\ov_{k}\! +\! \ov_{j'}v_{k})dt  \\
\ns\ds \q + 2\sum_{j,k,j',k'=1}^n
b^{jk}b^{j'k'}\ell_j (v_{j'}\ov_{kk'} +
\ov_{j'}v_{kk'})dt  \\
\ns\ds = - 2\!\! \sum_{j,k,j',k'=1}^n \!\[
b^{jk}b^{j'k'}\ell_j (v_{j'}\ov_{k} \!+\!
\ov_{j'}v_{k}) \]_{k'} dt \!+\!
2\!\!\!\sum_{j,k,j',k'=1}^n\!\!
b^{j'k'}(b^{jk}\ell_j)_{k'} (v_{j'}\ov_{k}\! + \!\ov_{j'}v_{k})dt  \\
\ns\ds \q + \!\!\!\sum_{j,k,j',k'=1}^n \!\[
b^{jk}b^{j'k'}\ell_j (v_{j'}\ov_{k'} \!+\!
\ov_{j'}v_{k'}) \]_k dt -
\!\!\sum_{j,k,j',k'=1}^n
(b^{jk}b^{j'k'}\ell_j)_k
(v_{j'}\ov_{k'} + \ov_{j'}v_{k'})dt \\
\ns\ds = - 2\!\! \sum_{j,k,j',k'=1}^n \!\[
b^{jk'}b^{j'k}\ell_j (v_{j'}\ov_{k'}\! +\!
\ov_{j'}v_{k'}) \]_{k} dt \!+\!
2\!\!\sum_{j,k,j',k'=1}^n\!\!
b^{jk'}(b^{j'k}\ell_{j'})_{k'} (v_{j}\ov_{k} \!+ \!\ov_{j}v_{k})dt  \\
\ns\ds \q + \!\!\sum_{j,k,j',k'=1}^n \[
b^{jk}b^{j'k'}\ell_j (v_{j'}\ov_{k'}\! +
\ov_{j'}v_{k'}) \]_k dt -
\!\!\sum_{j,k,j',k'=1}^n
(b^{jk}b^{j'k'}\ell_{j'})_{k'} (v_{j}\ov_{k} +
\ov_{j}v_{k})dt.
\end{array}
\end{equation}

Further, it holds that
\begin{equation}\label{th1eq10}
\begin{array}{ll}\ds
 \q I_1^3\overline{I_2^2}  +  I_2^2 \overline{I_1^3} \\
 \ns\ds = \Psi v  \overline{ \sum_{j,k=1}^n(b^{jk}v_j)_k }dt + \sum_{j,k=1}^n(b^{jk}v_j)_k \overline{ \Psi v }dt \\
 \ns\ds =
\sum_{j,k=1}^n
\[ \Psi b^{jk}(v_j \ov + \ov_j v) \]_k dt -
\sum_{j,k=1}^n
 \Psi b^{jk}(v_j \ov_k + \ov_j v_k)dt  \\
\ns\ds\q -\sum_{j,k=1}^n \Psi_k  b^{jk} ( v_j
\bar v + \bar v_j v) dt
 \\
 \ns\ds =
\sum_{j,k=1}^n
\[ \Psi b^{jk}(v_j \ov + \ov_j v) \]_k dt -
\sum_{j,k=1}^n
 \Psi b^{jk}(v_j \ov_k + \ov_j v_k)dt  \\
\ns\ds\q -\sum_{j,k=1}^n \[ b^{jk}\Psi_k |v|^2
\]_j dt + \sum_{j,k=1}^n (b^{jk}\Psi_k)_j
|v|^2dt.
\end{array}
\end{equation}

Finally, we have that
\begin{equation}\label{th1eq11}
\begin{array}{ll}\ds
\q I_1\overline{I_2^3}  +  I_2^3 \overline{I_1} \\
\ns\ds = -i\b \ell_t v \overline{Av}dt +
Av\overline{ (-i\b \ell_t v) }dt
 \\
\ns\ds =  - 2\sum_{j,k=1}^n (b^{jk}\ell_j
A|v|^2)_k dt + 2\[ \sum_{j,k=1}^n (b^{jk}\ell_j
A)_k + A\Psi \]|v|^2dt .
\end{array}
\end{equation}

{\bf Step 3.} Combining
(\ref{th1eq2})--(\ref{th1eq11}), we conclude
the desired formula (\ref{c2a1}).

\section{Carleman Estimate for Stochastic Schr\"{o}dinger Equations}

This section is devoted to the proof of Theorem
\ref{thcarleman est}.

\vspace{0.1cm}
{\em Proof of Theorem
{\ref{thcarleman est}}}: The proof is divided
into  three steps.

\medskip

\textbf{Step 1.} We choose $\b = 1$ and
$(b^{jk})_{1\leq j,k\leq n}$ to be the identity
matrix. Put $$ \d^{jk} =
\left\{\begin{array}{ll}\ds 1,&\mbox{ if }
j=k,\\ \ns\ds 0,&\mbox{ if } j\neq
k.\end{array}\right. $$   Applying Theorem
\ref{identity1} to the equation \eqref{system1}
with $\theta$ given by \eqref{lvarphi}, $z$
replaced by $y$ and $v = \theta z$. We obtain
that
\begin{equation}\label{identity2.1}
\begin{array}{ll}
\ds\q\theta\cP y  {\( i\b \ell_t \bar{v} -
2\sum_{j,k=1}^n b^{jk}\ell_j \bar{v}_k + \Psi
\bar{v}\)} + \theta\overline{\cP y} {\(- i\b
\ell_t v - 2\sum_{j,k=1}^n b^{jk}\ell_j v_k  +
\Psi v\)}\\
\ns \ds \q + \;dM + \div V  \\
\ns\ds  = 2\Big|- i\b \ell_t v -
2\sum_{j,k=1}^n b^{jk}\ell_j v_k  + \Psi
v\Big|^2dt +
\sum_{j,k=1}^nc^{jk}(v_k\ov_j+\ov_k
v_j) dt + D|v|^2dt  \\
\ns \ds \q + 2i\sum_{j=1}^n (\ell_{jt} +
\ell_{tj})(\ov_j v - v_j\ov)dt +
i(\Psi + \D \ell)(\ov dv - v d\ov)  \\
\ns \ds  \q + \ell_t dv d\ov + i\sum_{j=1}^n
\ell_j (d\ov_j dv - dv_j d\ov).
\end{array}
\end{equation}
Here
\begin{equation}\label{Id2eq1.1}
\begin{array}{ll}
M \3n& \ds = \b^2\ell_t |v|^2 +
i\b\sum_{j,k=1}^nb^{jk}\ell_j(\ov_kv-v_k\ov)\\
\ns & \ds = \ell_t |v|^2 + i\sum_{j=1}^n \ell_j
(\ov_j v - v_j \ov);
\end{array}
\end{equation}
\begin{equation}\label{Id2eq1.2}
\begin{array}{ll}
A \3n& \ds =\sum_{j,k=1}^nb^{jk}\ell_j\ell_k -
\sum_{j,k=1}^n(b^{jk}\ell_j)_k -\Psi \\
\ns & \ds  = \sum_{j=1}^n (\ell_j^2 -
\ell_{jj}) -\Psi;
\end{array}
\end{equation}
\begin{equation}\label{Id2eq1.3}
\begin{array}{ll}
D \3n & \ds =(\b^2\ell_t)_t
+\sum_{j,k=1}^n(b^{jk}\Psi_k)_j
+ 2\[\sum_{j,k=1}^n(b^{jk}\ell_j A)_k + A\Psi\]\\
\ns & \ds = \ell_{tt} + \sum_{j=1}^n \Psi_{jj}
+ 2\sum_{j=1}^n (\ell_j A)_j + 2 A\Psi;
\end{array}
\end{equation}
\begin{equation}\label{Id2eq1.4}
\begin{array}{ll}
c^{jk} \3n &\ds  =
\sum_{j',k'=1}^n\[2(b^{j'k}\ell_{j'})_{k'}b^{jk'}
- (b^{jk}b^{j'k'}\ell_{j'})_{k'}\Psi\] - b^{jk}\\
\ns & \ds  =
\[2(b^{kk}\ell_{k})_{j}b^{jj} -
\sum_{j'=1}^n (b^{jk}b^{j'j'}\ell_{j'})_{j'} -
b^{jk}\Psi\]\\
\ns & \ds = 2\ell_{jk} - \d^{jk}\D \ell -
\d^{jk}\Psi;
\end{array}
\end{equation}
and
\begin{equation}\label{Id2eq1.5}
\begin{array}{ll}
V_k \3n &\ds = -i
\b\sum_{j=1}^n\[b^{jk}\ell_j(vd\ov -\ov
 dv ) + b^{jk}\ell_t(v_j\ov-\ov_jv) dt\]\\
\ns & \ds \q - \Psi\sum_{j=1}^n b^{jk}(v_j\ov+\ov_jv) dt + \sum_{j=1}^n b^{jk}(2A\ell_j+\Psi_j)|v|^2 dt \\
\ns & \ds \q
+\sum_{j,j',k'=1}^n\(2b^{jk'}b^{j'k}-b^{jk}b^{j'k'}\)\ell_j(v_{j'}\ov_{k'}+\ov_{j'}v_{k'})
dt\\
\ns & \ds = -i\big[ \ell_k(vd\ov - \ov
dv) + \ell_t(v_j\ov -\ov_j v)dt \big] - \Psi(v_k\ov + \ov_k v)dt + (2A\ell_k + \Psi_k)|v|^2dt\\
\ns &  \ds \q + 2\sum_{j=1}^n \ell_j (\ov_j v_k
+ v_j \ov_k)dt - 2\sum_{j'=1}^n
\ell_k(v_j\ov_j)dt.
\end{array}
\end{equation}

\textbf{Step 2.} In this step, we estimate the
terms in the right-hand side of the equality
\eqref{identity2.1} one by one.

First, from the definition of $\ell$, $\f$(see
\eqref{lvarphi}) and the choice of $\psi$(see
\eqref{psi}), we have that
\begin{equation}\label{lt1}
\begin{array}{ll}\ds
|\ell_t| & \ds = \Big| s\frac{2(2t-T)}{t^3(T-t)^3}\big( e^{4\l\psi} - e^{5\l |\psi|_{L^\infty(G)}} \big)  \Big| \\
\ns& \ds \leq \Big| s\frac{2(2t-T)}{t^3(T-t)^3} e^{5\l |\psi|_{L^\infty(G)}} \Big|  \\
\ns &\ds \leq  \Big| s\frac{C}{t^3(T-t)^3} e^{5\l \psi} \Big|\\
\ns & \ds  \leq Cs\varphi^{1+\frac{1}{2}},
\end{array}
\end{equation}
and that
\begin{equation}\label{ltt1}
\begin{array}{ll}
\ds |\ell_{tt}| & \ds  = \Big|  s\frac{20t^2 - 20tT + 6T^2}{t^4(T-t)^4} \big( e^{4\l\psi} - e^{5\l |\psi|_{L^\infty(G)}} \big) \Big| \\
\ns& \ds \leq \Big|  s\frac{C}{t^4(T-t)^4}  e^{5\l |\psi|_{L^\infty(G)}}  \Big| \\
\ns& \ds \leq \Big|  s\frac{C}{t^4(T-t)^4}  e^{8\l  \psi }  \Big|\\
\ns  &\ds  \leq Cs\f^2\leq Cs\f^3.
\end{array}
\end{equation}

We choose  below $\Psi = -\D \ell$, then we
have that
\begin{eqnarray}\label{Id2eq2}
 A = \sum_{j=1}^m \ell_j^2 =  \sum_{j=1}^m \big(4s\l\f \psi )^2 =16s^2\l^2\varphi^2 |\nabla\psi|^2.
\end{eqnarray}
Hence, we find
\begin{equation}\label{B}
\begin{array}{ll}\ds
D \3n & \ds = \ell_{tt} + \sum_{j=1}^n
\Psi_{jj} + 2\sum_{j=1}^n (\ell_j
A)_j + 2 A\Psi \\
\ns & \ds  = \ell_{tt} + \D(\D\ell) + 2\sum_{j=1}^n\big(4s\l\f\psi_j 16s^2\l^2\f^2|\nabla\psi|^2\big)_j - 32s^2\l^2\f^2|\nabla\psi|^2\D \ell  \\
\ns & \ds  =  384s^3\l^4\varphi^3|\nabla\psi|^4
- \l^4\varphi O(s) - s^3\varphi^3 O(\l^3) +
\ell_{tt}.
\end{array}
\end{equation}
Recalling that $x_0\in (\mathbb{R}^n\setminus
\overline G)$, we know that
$$|\nabla\psi|>0\;\;\mbox{ in }\overline G.$$ From   \eqref{B} and  \eqref{ltt1}, we
know that there exists a $\l_0>0$ such that for
all $\l>\l_0$, one can find a constant $s_0 =
s_0(\l_0)$ so that for any $s>s_0$, it holds
that
\begin{equation}\label{B1}
D|v|^2 \geq
s^3\l^4\varphi^3|\nabla\psi|^4|v|^2.
\end{equation}
Since
$$
\begin{array}{ll}\ds
 c^{jk} = 2\ell_{jk} - \d^{jk}\D \ell - \d^{jk}\Psi  \\
\ns\ds\q\,\,\, = 32s\l^2\varphi\psi_j \psi_k +
16s\l\varphi\psi_{jk},
\end{array}
$$
we see that
\begin{equation}\label{cjk}
\begin{array}{ll}\ds
\q \ds\sum_{j,k=1}^n c^{jk}(v_j\ov_k + v_k\ov_j)\\
\ns \ds  = 32s\l^2\varphi\sum_{j,k=1}^n\psi_j
\psi_k(v_j\ov_k + v_k\ov_j) + 16s\l\varphi
\sum_{j,k=1}^n \psi_{jk}(v_j\ov_k + v_k\ov_j)\\
\ns\ds =
32s\l^2\varphi\[\sum_{j=1}^n(\psi_jv_j)\sum_{k=1}^n
(\psi_k \ov_k) +
\sum_{k=1}^n(\psi_kv_k)\sum_{j=1}^n (\psi_j
\ov_j)  \] + 32s\l\varphi \sum_{j=1}^n(v_j\ov_j
+ \ov_j v_j)\\
\ns   \ds  = 64s\l^2\varphi |\nabla\psi\cd\nabla v|^2 + 64 s\l\f |\nabla v|^2\\
\ns   \ds \geq  64 s\l\f |\nabla v|^2.
\end{array}
\end{equation}

Now we estimate the other terms in the
right-hand side of the equality
\eqref{identity2.1}. The first one satisfies
that
\begin{eqnarray}\label{ltj}
\begin{array}
{ll} \ds  2i\sum_{j=1}^n (\ell_{jt} +
\ell_{tj})(\ov_j v - v_j\ov) & \ds =
4i\sum_{j=1}^n s\l\psi_j \ell_t(\ov_j v - \ov v_j)\\
\ns & \ds \leq 2 s\varphi |\nabla v|^2 + 2
s\l^2\varphi^3 |\nabla\psi|^2|v^2|.
\end{array}
\end{eqnarray}

The second one reads
\begin{eqnarray}\label{liiPsi}
i(\Psi + \D \ell)(\ov dv - v d\ov) = 0.
\end{eqnarray}

For the estimate of the third and the fourth
one, we need to take mean value and get that
\begin{equation}\label{dvdov}
\begin{array}{ll}\ds
 \mathbb{E}\big(\ell_t dv d\ov\big) \3n& \ds= \mathbb{E}\big[\ell_t(\theta \ell_t ydt + \theta
dy)(\overline{\theta \ell_t ydt + \theta
dy)}\big] = \mathbb{E}(\ell_t \theta^2 dy
d\bar{y})
\\
\ns & \ds \leq 2s\theta^2
\varphi^{\frac{3}{2}}\mathbb{E}( a_3^2|y|^2 +
g^2)dt.
\end{array}
\end{equation}
Here we utilize inequality \eqref{lt1}.

\vspace{0.1cm}

Since
$$
\begin{array}{ll}\ds
\mathbb{E}(d\ov_j dv) & = \mathbb{E}\big[\overline{\big( \theta \ell_t v dt + \theta dy  \big)}_j \big( \theta \ell_t v dt + \theta dy  \big)\big] \\
\ns& \ds  = \mathbb{E} \big[\, \overline{(\theta dy)}_j (\theta dy) \big]\\
\ns& \ds =  \mathbb{E} \big[\, \overline{\big( s\l\f\psi_j\theta dy + \theta dy_j   \big)}\theta dy  \big]\\
\ns & \ds  = s\l\f\psi_j\theta^2 \mathbb{E}d\bar y dy + \theta^2 \mathbb{E}d\bar y_j dy \\
\ns & \ds  =  s\l\f\psi_j\theta^2
\mathbb{E}|a_3y + g|^2dt + \theta^2
\mathbb{E}\big[\,\overline{ (a_3 y + g) }_j
(a_3 y + g) \big]dt
\end{array}
$$
and
$$
\begin{array}{ll}\ds
\q\theta^2 \mathbb{E}\big[\,\overline{ (a_3 y +
g) }_j
(a_3 y + g) \big]dt\\
\ns\ds = \theta^2 \mathbb{E}\big[(\overline{a_3
y})_j (a_3 y) + (\overline{a_3 y})_j g + (a_3 y
)\bar g_j + g\bar g_j \big] dt\\
\ns\ds =\theta^2 \mathbb{E}\big[(\overline{a_3
y})_j (a_3 y) + (\overline{a_3 y})_j g + g\bar
g_j \big] dt + [\mE\theta^2(a_3 y )\bar g]_j \\
\ns\ds \q - s\l\f\psi_j\theta^2\mE(a_3 y \bar
g)-\th^2\mE[(a_3y)_j]\bar g,
\end{array}
$$
we get that
$$
\begin{array}{ll}\ds
\mathbb{E}(d\ov_j dv) \3n&\ds=
s\l\f\psi_j\theta^2 \mathbb{E}|a_3y + g|^2dt +
\theta^2 \mathbb{E}\big[(\overline{a_3 y})_j
(a_3 y) + (\overline{a_3 y})_j g   + g\bar g_j
\big]
dt\\
\ns&\ds \q + \mE(\theta^2 a_3 y \bar g)_j -
s\l\f\psi_j\theta^2\mE(a_3 y \bar
g)-\th^2\mE[(a_3y)_j\bar g].
\end{array}
$$
Similarly, we can get that
$$
\begin{array}{ll}\ds
\mathbb{E}(dv_j d\ov)\3n& \ds=
s\l\f\psi_j\theta^2 \mathbb{E}|a_3y + g|^2dt +
\theta^2 \mathbb{E}\big[(\overline{a_3 y}) (a_3
y)_j  + (a_3 y )_j\bar g + g_j\bar g \big]
dt\\
\ns&\ds \q + \mE(\theta^2 \overline{a_3 y} g)_j
- s\l\f\psi_j\theta^2\mE(\overline{a_3 y}
 g)-\th^2\mE[(\overline{a_3 y})_j g].
\end{array}
$$
Therefore, fourth one enjoys that
\begin{equation}\label{dvjdv}
\begin{array}{ll}
\ds \q i\mathbb{E}\sum_{j=1}^n \ell_j (d\ov_j
dv -
dv_j d\ov)\\
\ns\ds =   s\l\varphi\sum_{j=1}^n \psi_j
\[\mathbb{E}\big(d\ov_j dv\big) - \mathbb{E}\big(dv_j d\ov\big) \]  \\
\ns \ds =   s\l\varphi \psi \sum_{j=1}^n \psi_j
\theta^2 \mathbb{E}\Big\{\big[(\overline{a_3
y})_j (a_3 y) + (\overline{a_3 y})_j g  + g\bar
g_j  - s\l\f\psi_j a_3 y \bar
g -  (a_3y)_j \bar g\big]\\
\ns \ds \q - \big[(\overline{a_3 y}) (a_3 y)_j
+ (a_3 y )_j\bar g + g_j\bar g  - s\l\f\psi_j
(\overline{a_3 y} g)- [(\overline{a_3 y})_j
g]\big]\Big\} dt \\
\ns\ds \q + s\l\varphi \psi \sum_{j=1}^n \psi_j
\mathbb{E}\big(\th^2 a_3 y \bar g - \theta^2
\overline{a_3 y} g \big)_j.
\end{array}
\end{equation}

\textbf{Step 3.} Integrating the equality
\eqref{identity2.1} in $Q$,  taking mean value
in both sides, and noting
\eqref{Id2eq2}--\eqref{dvjdv},  we obtain that
\begin{equation}\label{inep1}
\begin{array}{ll}
\ds \q\mathbb{E}\int_Q \Big(s^3\l^4\varphi^3
|v|^2 \!+\! s\l^2\varphi |\nabla v|^2\Big) dxdt
+ 2\mathbb{E}\int_Q \Big|\!- i\b
\ell_t v - 2\!\!\sum_{j,k=1}^n  b^{jk}\ell_j v_k  + \!\Psi v\Big|^2dxdt\\
\ns \ds  \leq  \mathbb{E}\int_Q \Big\{
\theta\cP y {\Big( i\b \ell_t \bar{v}\! -\!
2\!\! \sum_{j,k=1}^n\! b^{jk}\ell_j \bar{v}_k
+ \!\Psi \bar{v}\Big)} + \theta\overline{\cP y}
{\Big(\!-\! i\b \ell_t v\! - \!2\!
\sum_{j,k=1}^n\!
b^{jk}\ell_j v_k + \!\Psi v\Big)} \Big\}dx\\
\ns \ds \q  +\; C\mathbb{E}\int_Q
\theta^2\Big[s^2\l^2 \varphi^2(a_3^2|y|^2 +
g^2) + a_3^2|\nabla y|^2 + |\nabla a_3|^2 y^2
+ |\nabla g|^2\Big]  dxdt\\
\ns\ds \q + \;\mathbb{E}\int_Q dM dx +
\mathbb{E}\int_Q \div V dx.
\end{array}
\end{equation}

Now we analyze the  terms in the right-hand
side of the inequality \eqref{inep1} one by
one.

The first term satisfies that
\begin{equation}\label{intprin}
\begin{array}{ll} \ds \mathbb{E}\int_Q \Big\{ \theta\cP y
{\Big( i\b \ell_t \bar{v} - 2\sum_{j,k=1}^n
b^{jk}\ell_j \bar{v}_k  + \Psi
\bar{v}\Big)}\\
\ns  \ds \q +\; \theta\overline{\cP y} {\Big(-
i\b \ell_t v -
2\sum_{j,k=1}^n b^{jk}\ell_j v_k + \Psi v\Big)} \Big\}dx \\
\ns  \ds  =  \ds \mathbb{E}\int_Q \Big\{ \theta
(a_1 \cdot \nabla y + a_2 y + f) {\( i\b \ell_t
\bar{v} - 2\sum_{j,k=1}^n b^{jk}\ell_j
\bar{v}_k + \Psi \bar{v}\)}\\
\ns\ds \q +\; \theta {(a_1 \cdot \nabla \bar{y}
+ \overline{a_2 y} + \bar{f})} {\Big(- i\b
\ell_t v - 2\sum_{j,k=1}^n b^{jk}\ell_j v_k  +
\Psi
v\Big)} \Big\}dxdt\nonumber\\
\ns \ds \leq  2\mathbb{E}\int_Q
\Big\{\theta^2\big|a_1 \cdot \nabla y + a_2 y +
f\big|^2 + \Big|- i\b \ell_t v -
2\sum_{j,k=1}^n b^{jk}\ell_j v_k + \Psi
v\Big|^2 \Big\}dxdt.
\end{array}
\end{equation}

From the definition of $\theta$, we know that
$v(0)=v(T)=0$. Hence, it holds that
\begin{equation}\label{idm}
\int_Q dM dx = 0.
\end{equation}
By means of Stokes' Theorem, we have that
\begin{eqnarray}\label{intV}
\begin{array}
{ll} \ds \mathbb{E}\int_Q \div V dx  \3n&\ds =
\ds \mathbb{E}\int_{\Si}
2\sum_{k=1}^n\sum_{j=1}^n\Big[
 \ell_j\big(\ov_j v_k + v_j \ov_k\big)\nu^k - \ell_k \nu_k v_j \ov_j
\Big]d\Si\\
\ns  &\ds =  \ds \mathbb{E}\int_{\Si}
\Big(4\sum_{j=1}^n \ell_j \nu_j \Big| \frac{\pa
v}{\pa \nu} \Big|^2 - 2\sum_{k=1}^n \ell_k
\nu_k \Big|
\frac{\pa v}{\pa \nu} \Big|^2\Big) d\Si\\
\ns  &= \ds  \mathbb{E}\int_{\Si} 2\sum_{k=1}^n
\ell_k \nu_k \Big|
\frac{\pa v}{\pa \nu} \Big|^2 d\Si \\
\ns &\ds \leq  C\mathbb{E}\int_0^T \int_{\G_0}
\theta^2 s\l\varphi \Big| \frac{\pa y}{\pa \nu}
\Big|^2 d\G dt.
\end{array}
\end{eqnarray}

By (\ref{inep1})--(\ref{intV}), we have that
\begin{eqnarray}\label{car1}
\begin{array}{ll}
\q \ds \mathbb{E}\int_Q \Big(s^3\l^4\varphi^3
|v|^2 + s\l\varphi
|\nabla v|^2\Big) dxdt \\
\ns   \ds \leq C\,\mathbb{E}\int_Q \theta^2
|a_1 \cdot \nabla y + a_2 y + f|^2 dxdt +
C\,\mathbb{E}\int_0^T\int_{\G_0}\theta^2
s\l\varphi\Big| \frac{\pa y}{\pa \nu}\Big|^2d\G dt\\
\ns \ds   \q  +\, C\mathbb{E}\int_Q
\theta^2\Big[s^2\l^2 \varphi^2\big(a_3^2|y|^2 +
g^2\big) + a_3^2|\nabla y|^2 + |\nabla a_3|^2
y^2 + |\nabla g|^2\Big] dxd t.
\end{array}
\end{eqnarray}

Noting that $y_i = \theta^{-1}(v_i - \ell_i v)
= \theta^{-1}(v_i - s\l\varphi\psi_i v)$, we
get
\begin{equation}\label{vtoy}
\theta^2\big(|\nabla y|^2 + s^2\l^2\varphi^2
|y|^2\big)\leq C\big(|\nabla v|^2 +
s^2\l^2\varphi^2 |v|^2\big).
\end{equation}

Therefore, it follows from (\ref{car1}) that
\begin{equation}\label{car2}
\begin{array}{ll}
\ds \q\mathbb{E}\int_Q \Big(s^3\l^4\varphi^3
|y|^2 + s\l\varphi
|\nabla y|^2\Big) dxdt \\
\ns \ds   \leq  C\mathbb{E}\int_Q \theta^2\Big(
|a_1|^2 || \nabla y|^2 + a_2^2 |y|^2 +
|f|^2\Big) dxdt +
C\mathbb{E}\int_0^T\int_{\G_0}\theta^2
s\l\varphi\Big| \frac{\pa
y}{\pa \nu}\Big|^2d\G dt \\
\ns  \ds \q + C\mathbb{E}\int_Q
\theta^2\Big[s^2\l^2 \varphi^2\big(a_3^2|y|^2 +
g^2\big) + a_3^2|\nabla y|^2 + |\nabla a_3|^2
y^2 + |\nabla g|^2\Big]  dxdt.
\end{array}
\end{equation}

Taking   $\l_1 =\l_0$ and $s_1 = \max(s_0,
Cr_1)$, and utilizing the inequality
\eqref{car2}, we conclude the desired
inequality \eqref{carleman est}.

On the other hand, if $g\in
L^2_\cF(0,T;H^1(G;\mathbb{R}))$, then $g\bar
g_j - g_j\bar g=0$ for $j=1,\cds,n$. Thus, from
\eqref{Id2eq2}--\eqref{dvjdv}, we get
\begin{equation}\label{inep1z}
\begin{array}{ll}
\ds \q\mathbb{E}\int_Q \Big(s^3\l^4\varphi^3
|v|^2 + s\l^2\varphi |\nabla v|^2\Big) dxdt +
2\mathbb{E}\int_Q \Big|\!- i\b
\ell_t v \!- \!2\!\sum_{j,k=1}^n b^{jk}\ell_j v_k \! +\! \Psi v\Big|^2dxdt\\
\ns \ds  \leq  \mathbb{E}\int_Q \Big\{
\theta\cP y {\Big( i\b \ell_t \bar{v} \!-\!
2\!\!\sum_{j,k=1}^n b^{jk}\ell_j \bar{v}_k  +
\Psi \bar{v}\Big)} + \theta\overline{\cP y}
{\Big(\!\!-\! i\b \ell_t v\! -\!
2\!\!\sum_{j,k=1}^n
b^{jk}\ell_j v_k\! + \!\Psi v\Big)} \Big\}dx\\
\ns \ds \q  +\; C\mathbb{E}\int_Q
\theta^2\Big[s^2\l^2 \varphi^2\big(a_3^2|y|^2 +
g^2\big) + a_3^2|\nabla y|^2  + |\nabla a_3|^2
y^2
\Big]  dxdt +  \mathbb{E}\int_Q dM dx \\
\ns\ds \qq  + \mathbb{E}\int_Q \div V dx.
\end{array}
\end{equation}
Then, by a similar argument, we find that
\begin{equation}\label{car2z}
\begin{array}{ll}
\ds \q\mathbb{E}\int_Q \Big(s^3\l^4\varphi^3
|y|^2 + s\l\varphi
|\nabla y|^2\Big) dxdt \\
\ns \ds   \leq  C\mathbb{E}\int_Q \theta^2\Big(
|a_1|^2 || \nabla y|^2 + a_2^2 |y|^2 +
|f|^2\Big) dxdt +
C\mathbb{E}\int_0^T\int_{\G_0}\theta^2
s\l\varphi\Big| \frac{\pa
y}{\pa \nu}\Big|^2d\G dt \\
\ns  \ds \q + C\mathbb{E}\int_Q
\theta^2\Big[s^2\l^2 \varphi^2\big(a_3^2|y|^2 +
g^2\big) + a_3^2|\nabla y|^2 + |\nabla a_3|^2
y^2 \Big]  dxdt.
\end{array}
\end{equation}
Now taking   $\l_1 =\l_0$ and $s_1 = \max(s_0,
Cr_1)$, and using the inequality \eqref{car2z},
we obtain the desired inequality
\eqref{carleman est1}.

\section{Proof of Theorem \ref{observability}}

In this section, we prove Theorems
\ref{observability}, by means of Theorem
\ref{thcarleman est}.

{\em Proof of Theorem \ref{observability}}: By
means of the definition of $\ell$ and
$\theta$(see \eqref{lvarphi}), it holds that
\begin{eqnarray}\label{final1}
\begin{array}{ll}
\ds \q\mathbb{E}\int_Q \theta^2\Big(\varphi^3
|y|^2 + \varphi
|\nabla y|^2\Big) dxdt\\
\ns \ds \geq
\min_{x\in\overline{G}}\Big(\varphi\Big(\frac{T}{2},x\Big)
\theta^2\Big(\frac{T}{4},x\Big)\Big)\mathbb{E}\int_{\frac{T}{4}}^{\frac{3T}{4}}\int_G\big(|y|^2+|\nabla
y|^2\big)dxdt,
\end{array}
\end{eqnarray}

\begin{equation}\label{final2}
\begin{array}{ll} \ds
\q\mathbb{E}\int_Q \theta^2\big(|f|^2 +
\varphi^2|g|^2 + |\nabla
g|^2\big)dxdt \\
\ns \ds\leq  \max_{(x,t)\in
\overline{Q}}\big(\varphi^2(t,x)\theta^2(t,x)\big)\mathbb{E}\int_Q\big(|f|^2
+ |g|^2 + |\nabla g|^2\big)dxdt
\end{array}
\end{equation}
and that
\begin{equation}\label{final3}
\mathbb{E}\int_0^T\int_{\G_0}\theta^2
\varphi\Big| \frac{\pa y}{\pa \nu}\Big|^2d\G dt
\leq \max_{(x,t)\in
\overline{Q}}\big(\varphi(t,x)\theta^2(t,x)\big)\mathbb{E}\int_0^T\int_{\G_0}
\Big| \frac{\pa y}{\pa \nu}\Big|^2d\G dt.
\end{equation}

From  \eqref{carleman est} and
\eqref{final1}--(\ref{final3}), we deduce that
\begin{equation}\label{final4}
\begin{array}{ll} \ds
\q\mathbb{E}\int_{\frac{T}{4}}^{\frac{3T}{4}}\int_G\big(|y|^2+|\nabla
y|^2\big)dxdt\\
\ns \ds\leq   C r_1 \frac{\max_{(x,t)\in
\overline{Q}}\Big(\varphi^2(t,x)\theta^2(t,x)\Big)}{\min_{x\in\overline{G}}\Big(\varphi(\frac{T}{2},x)\theta^2(\frac{T}{4},x)\Big)}\\
\ns \ds \q\times\left\{
\mathbb{E}\int_Q\big(|f|^2 + |g|^2 + |\nabla
g|^2\big)dxdt + \mathbb{E}\int_0^T\int_{\G_0}
\Big| \frac{\pa y}{\pa
\nu}\Big|^2d\G dt\right\}\\
\ns  \ds \leq  e^{ Cr_1 }\left\{
\mathbb{E}\int_Q\big(|f|^2 + |g|^2 + |\nabla
g|^2\big)dxdt + \mathbb{E}\int_0^T\int_{\G_0}
\Big| \frac{\pa y}{\pa \nu}\Big|^2d\G
dt\right\}.
\end{array}
\end{equation}

Utilizing (\ref{final4}) and
(\ref{energyesti1}), we obtain that
\begin{equation}\label{final5}
\begin{array}{ll}  \ds
\q\mathbb{E}\int_G\big(|y_0|^2 + |\nabla y_0|^2\big)dx \\
\ns \ds \leq   e^{C r_1 }\left\{
\mathbb{E}\int_Q\big(|f|^2 + |\nabla f|^2 +
|g|^2 + |\nabla g|^2\big)dxdt +
\mathbb{E}\int_0^T\int_{\G_0} \Big| \frac{\pa
y}{\pa \nu}\Big|^2d\G dt\right\},
\end{array}
\end{equation}
which concludes Theorem \ref{observability}
immediately.\endproof

\medskip

\section{Two applications}\label{Sec app}

This section is addressed to applications of
the observability/Carleman estimates shown in
Theorems \ref{observability}--\ref{thcarleman
est}.

We first study a state observation problem for
semilinear stochastic Schr\"{o}dinger
equations. Let us consider the following
equation:
\begin{equation}\label{system2}
\!\!\left\{
\begin{array}{ll}
\ds idz + \D zdt=\big[\,a_1 \cdot \nabla z +
a_2 z +F_1(|z|)\big]dt + \big[a_3 z
+F_2(|z|)\big]
dB(t)&\mbox{ in } Q,\\
\ns\ds
z=0&\mbox{ on }\Si,\\
\ns\ds z(0)=z_0, &\mbox{ in }G.
\end{array}
\right.
\end{equation}
Here $a_i$ ($i=1,2,3$) are given as in
\eqref{coai}, $F_1(\cd)\in C^1(\mathbb{R};
\mathbb{C})$ with $F(0)=0$ and $F_2(\cd)\in
C^1(\mathbb{R}; \mathbb{R})$ are two known
nonlinear global Lipschitz continuous functions
with Lipschitzian constant $L$, while the
initial datum $z_0\in L^2(\O,\cF_0,P;
H_0^1(G))$ is unknown. The solution to the
equation \eqref{system2} is understood similar
to Definition \ref{def1}.

\begin{remark}
From the choice of $F_1$ and $F_2$, one can
easily show that the equation \eqref{system2}
admits a unique solution $z\in H_T$ by the
standard fixed point argument. We omit the
proof here.
\end{remark}

The state observation problem associated to the
equation \eqref{system2} is as follows.

\vspace{0.1cm}
\begin{itemize}

\item {\bf Identifiability}. Is the solution
$z\in H_T$ (to \eqref{system2}) determined
uniquely by the observation $\ds\frac{\pa
z}{\pa\nu}\Big|_{(0,T)\times \G_0}$?

\vspace{0.1cm}

\item {\bf Stability}. Assume that two
solutions $z$ and $\hat z$ (to \eqref{system2})
are given, and let $\ds\frac{\pa
z}{\pa\nu}\Big|_{(0,T)\times \G_0}$ and
$\ds\frac{\pa \hat z}{\pa\nu}\Big|_{(0,T)\times
\G_0}$ be the corresponding observations. Can
we find a positive constant $C$ such that
$$
|\!| z-\hat z |\!| \leq C\|\!\| \frac{\pa
z}{\pa\nu}-\frac{\pa \hat z}{\pa\nu} \|\!\|,
$$
with appropriate norms in both sides?

\vspace{0.1cm}

\item {\bf Reconstruction}. Is it possible to
reconstruct $z\in H_T$ to \eqref{system2}, in
some sense,  from the observation $\ds\frac{\pa
z}{\pa\nu}\Big|_{(0,T)\times \G_0}$?

\end{itemize}

The state observation problem for systems
governed by deterministic partial differential
equations is studied extensively (See
\cite{Kli,Li1,Yamamoto} and the rich references
therein). However, the stochastic case attracts
very little attention. To our best knowledge,
\cite{Zhangxu3} is the only published paper
addressing this topic. In that paper, the
author studied the state observation problem
for semilinear stochastic wave equations. By
means of Theorem \ref{observability}, we can
give positive answers to the above first and
second questions.

We claim that $\frac{\pa
z}{\pa\nu}|_{(0,T)\times \G_0}\in
L^2_\cF(0,T;L^2(\G_0))$ (and therefore, the
observation makes sense). Indeed, from the
choice of $F_1$,  it follows that
$$
\begin{array}{ll}\ds
\mathbb{E}\int_0^T\!\int_G \big|\n
\big(F_1(|z|)\big)\big|^2dxdt \3n& \ds=
\mathbb{E}\int_0^T\!\int_G \| F_1' (|z|)\n |z|
\|^2dxdt \leq L\mathbb{E}\int_0^T\!\int_G
\big|\n |z|\big|^2dxdt\\
\ns&\ds \leq L\mathbb{E}\int_0^T\int_G \big|\n
 z\big|^2dxdt,
\end{array}
$$
and
$$
F(|z(t,\cd)|)=0 \q\mbox{ on } \G \mbox{ for
a.e. }\  t\in [0,T].
$$
Hence,
$$
F_1(|z|)\in L^2_{\cF}(0,T;H_0^1(G)) \mbox{ for
any } z\in H_{T}.
$$
Similarly,
$$
F_2(|z|)\in L^2_{\cF}(0,T;H^1(G)) \mbox{ for
any } z\in H_{T}.
$$
Consequently, by Proposition~\ref{hregularity},
we find that $\frac{\pa
z}{\pa\nu}|_{(0,T)\times\G_0}\in
L^2_\cF(0,T;L^2(\G_0))$.

Now, we define a nonlinear map as follows:
$$
\left\{
\begin{array}{ll}\ds \cM:\
L^2(\O,\cF_0,P;H_0^1(G))\to
L^2_{\cF}(0,T;L^2(\G_0)),\\
\ns\ds \cM(z_0)= {\frac{\pa z}{\pa
\nu}}\Big|_{(0,T)\times \G_0},
\end{array}
\right.
$$
where $z$ solves the equation \eqref{system2}.
We have the following result.

\begin{theorem}\label{th2}
There exists a constant $\wt C=\wt C(L,T,G)>0$
such that for any $z_0, \hat z_0\in
L^2(\O,\cF_0,P;H_0^1(G))$, it holds that
\begin{equation}\label{th2eq1}
|z_0-\hat z_0|_{L^2(\O,\cF_0,P;L^2(G))} \le \wt
C|\cM(z_0) -\cM(\hat
z_0)|_{L^2_{\cF}(0,T;L^2(\G_0))},
\end{equation}
where $\hat z=\hat z(\cd\, ;\hat z_0)\in H_{T}$
is the solution to \eqref{system2} with $z_0$
replaced by $\hat z_0$.
\end{theorem}

\begin{remark}  From the
well-posedness of the equation \eqref{system2},
Theorem~\ref{th2} indicates that the state
$z(t)$ of \eqref{system2} (for $t\in [0,T]$)
can be uniquely determined from the observed
boundary data $\ds{\frac{\pa z}{\pa \nu}}
\Big|_{(0,T)\t\G_0}$, $P$- a.s., and
continuously depends on it. Therefore, we
answer the first and second questions for state
observation problem for the system
\eqref{system2} positively.
\end{remark}

{\it Proof of Theorem~\ref{th2}}\,: Set
$$
y=z-\hat z.
$$
Then, it is easy to see that $y$ satisfies
$$
\left\{
\begin{array}{ll}\ds
idy + \D y dt  =  \big[ a_1 \cdot \nabla y +
a_2 y +F_1(|z|)-F_1(|\hat z|) \big]dt \\
\ns\ds \hspace{2.2cm} + \big[ a_3 y +
F_2(|z|)-F_2(|z|) \big]dB(t) &\mbox{ in } Q,\\
\ns\ds y=0 &\mbox{ on }\Si,\\
\ns\ds y(0)=z_0-\hat z_0 &\mbox{ in } G.
\end{array}
\right.
$$
Also, it is clear that
$$
F_1(|z|)-F_1(|\hat z|)\in
L^2_{\cF}(0,T;H_0^1(G))
$$
and
$$
F_2(|z|)-F_2(|\hat z|)\in
L^2_{\cF}(0,T;H^1(G)).
$$
Hence, we know that $y$ solves the equation
\eqref{system1} with
$$
\left\{
\begin{array}{ll}\ds f=F_1(|z|)-F_1(|\hat z|),\\
\ns\ds g=F_2(|z|)-F_2(|\hat z|).
\end{array}
\right.
$$

By means of the inequality \eqref{carleman
est1} in Theorem \ref{thcarleman est},  there
exist an $s_1>0$ and a $\l_1>0$ so that for all
$s\geq s_1$ and $\l\geq \l_1$, it holds that
$$
\begin{array}{ll}\ds
\q\mathbb{E}\int_Q
\theta^2\Big(s^3\l^4\varphi^3 |y|^2 +
s\l\varphi
|\nabla y|^2\Big) dxdt \\
\ns \ds \leq  C \Big\{\mathbb{E}\int_Q \theta^2
\Big(|f|^2 +
 s^2\l^2\varphi^2 |g|^2  \Big)dxdt +
\mathbb{E}\int_0^T\int_{\G_0}\theta^2
s\l\varphi\Big| \frac{\pa y}{\pa \nu}\Big|^2d\G
dt \Big\}.
\end{array}
$$
By the choice of $f$, we see that
$$
\begin{array}{ll}\ds
\mathbb{E}\int_Q \theta^2 |f|^2dxdt \3n&\ds\leq
\mathbb{E}\int_Q \theta^2 |F_1(|z|)-F_1(|\hat
z|)|^2dxdt \leq L\mathbb{E}\int_Q \theta^2
(|z|-|\hat z|)^2dxdt\\
\ns&\ds \leq L\mathbb{E}\int_Q \theta^2 |z -
\hat z|^2dxdt \leq L\mathbb{E}\int_Q \theta^2
|y|^2dxdt.
\end{array}
$$
Similarly,
$$
s^2\l^2\mathbb{E}\int_Q \theta^2 \f^2 |g|^2dxdt
\leq  L\mathbb{E}\int_Q \theta^2\f^2 |y|^2dxdt.
$$
Hence, we obtain that
$$
\begin{array}{ll}\ds
\q\mathbb{E}\int_Q
\theta^2\Big(s^3\l^4\varphi^3 |y|^2 +
s\l\varphi
|\nabla y|^2\Big) dxdt \\
\ns \ds \leq  C\Big\{L \mathbb{E}\int_Q
\theta^2 \Big(|y|^2 +
 s^2\l^2\varphi^2 |y|^2  \Big)dxdt +
\mathbb{E}\int_0^T\int_{\G_0}\theta^2
s\l\varphi\Big| \frac{\pa y}{\pa \nu}\Big|^2d\G
dt \Big\}.
\end{array}
$$
Thus, there is a $\l_2\geq \max\{\l_1, CL\}$
such that for all $s\geq s_1$ and $\l\geq
\l_2$, it holds that
\begin{equation}\label{10.30eq1}
\mathbb{E}\int_Q \theta^2\Big(s^3\l^4\varphi^3
|y|^2 + s\l\varphi |\nabla y|^2\Big) dxdt \leq
C \mathbb{E}\int_0^T\int_{\G_0}\theta^2
s\l\varphi\Big| \frac{\pa y}{\pa \nu}\Big|^2d\G
dt.
\end{equation}

Further, similar to the proof of the inequality
\eqref{Eyt}, we can obtain that for any $0\leq
t \le s\leq T$, it holds
\begin{equation}\label{Eyt1}
\begin{array}{ll}\ds
\mathbb{E}| y(t)|^2_{ L^2(G)} - \mathbb{E}|
y(s)|^2_{ L^2(G)}\3n &\ds \leq
2\mathbb{E}\int_t^s\int_G \Big[ |f|^2 + |g|^2
\Big]dxd\si\\
\ns&\ds \leq CL\mathbb{E}\int_t^s\int_G |y|^2
dxd\si.
\end{array}
\end{equation}
Then, by Gronwall's inequality, we find that
\begin{equation}\label{Eyt2}
\mathbb{E}| y(t)|^2_{ L^2(G)} \leq
e^{CL}\mathbb{E} | y(s)|^2_{ L^2(G)}, \mbox{
for any } 0\leq t\leq s \leq T.
\end{equation}
Combining \eqref{10.30eq1} and \eqref{Eyt2},
similar as the derivation of the inequality
\eqref{final5}, we obtain the inequality
\eqref{th2eq1}.
\endproof

\vspace{0.2cm}

Now we consider the unique continuation
property for the equation \eqref{system1}.
There are numerous works on the unique
continuation property for deterministic partial
differential equations. The study in this
respect began at the very beginning of the 20th
century; while a climax appeared in the last
1950-70's. The most powerful tools at that
period is the local Carleman estimate (See
\cite{Hor1} for example). Nevertheless, most of
the existing works are devoted to the local
unique continuation property at that time. In
the recent 20 years, motivated by
Control/Inverse Problems of partial
differential equations, the study of the global
unique continuation is very active (See
\cite{Castro-Zuazua,Zhangxu1,Zhang-Zuazua} and
the rich references therein). Compared with the
fruitful works on the unique continuation
property in the deterministic settings, there
exist few results for stochastic partial
differential equations. As far as we know,
\cite{Zhangxu4,Zhangxu2} are the only two
published articles addressed to this topic, and
there is no result on the global unique
continuation property for stochastic
Schr\"{o}dinger equations in the previous
literature.

We remark that the powerful approach based on
local Carleman estimate in the deterministic
settings is very hard to apply to the
stochastic counterpart. Indeed, the usual
approach to employ local Carleman estimate for
the unique continuation needs to localize the
problem. Unfortunately, one cannot simply
localize the problem as usual in the stochastic
situation, since the usual localization
technique may change the adaptedness of
solutions, which is a key feature in the
stochastic setting. In this paper, as a
consequence of Theorem \ref{observability}
(which is based on the global Carleman estimate
established in Theorem \ref{thcarleman est}),
we obtain the following unique continuation
property for solutions to the equation
\eqref{system1}.

\vspace{0.1cm}

\begin{theorem}\label{ucp}
For any $\e>0$, let
$$O_\e([0,T]\t\G_0)\=\Big\{(x,t)\in Q
:\,\dist\big((x,t),[0,T]\t\G_0\big)\leq \e
\Big\}.$$ Let $f=g=0$, $P$-a.s. For any $y$
which solves the equation \eqref{system1}, if
 \bel{zx11}
 y = 0\ \ \hbox{ in }O_\e([0,T]\t\G_0), \ P\hbox{-a.s.},
 \ee
then $y=0$ in $Q$, $P$-a.s.
\end{theorem}

{\em Proof}\,: By \eqref{zx11}, we see that
$\ds\frac{\pa y}{\pa\nu}=0$ on $(0,T)\t\G_0$,
$P$-a.s. Hence, by means of Theorem
\ref{observability}, we find that $y(0)=0$ in
$L^2(\O,\cF_0,P;H_0^1(G))$. Consequently, we
conclude that $y=0$ in $Q$, $P$-a.s.
\endproof

\section{Further comments and open problems}

The subject of this paper is full of open
problems. Some of them seem to be particularly
relevant and could need important new ideas and
further developments:

\begin{itemize}

\item  {\bf Observability estimate for backward stochastic Schr\"{o}dinger equations}

Compared with Theorem \ref{observability}, it
is more interesting and difficult  to establish
the boundary observability estimate for
backward stochastic Schr\"{o}dinger equations.
More precisely, let us consider the following
backward stochastic Schr\"{o}dinger equation:
\begin{equation}\label{bsystem1}
\!\!\left\{
\begin{array}{lll}
\ds idu + \D u dt = (a_1\cd \nabla u + a_2 u + f)dt + (a_3 u + U+ g)dB(t) &\mbox{ in } Q,\\
\ns\ds u = 0 &\mbox{ on } \Si,\\
\ns\ds u(T) = u_T &\mbox{ in } G.
\end{array}
\right.
\end{equation}
Here the final state $u_T\in
L^2(\O,\cF_T,P;H_0^1(G))$ and
$\{\cF_t\}_{t\geq 0}$ is the natural filtration generated by $\{B(t)\}_{t\geq 0}$.  We expect the following result:\\

{\it Under the assumptions
\eqref{G0}--\eqref{cA}, any solution of the
equation \eqref{bsystem1}  satisfies that
\begin{equation} \label{bobser esti2}
\begin{array}{ll}\ds
 \q |u_T|_{L^2(\Omega,{ \mathcal{F}}_T, P; H_0^1(G))} \\
\ns\ds \leq  e^{C r_1} \Big(\Big|\frac{\partial
u}{\partial \nu}\Big |_{L^2_{ \mathcal{
F}}(0,T;L^2(\Gamma_0))} + |f|_{L^2_{ \mathcal{
F}}(0,T;H_0^1(G))} + |g|_{L^2_{ \mathcal{
F}}(0,T;H^1(G))}\Big),
\end{array}
\end{equation}
or at least,
\begin{equation} \label{bobser esti3}
\begin{array}{ll}\ds
 \q |u(0)|_{L^2(\Omega,{ \mathcal{F}}_0, P; H_0^1(G))} \\
\ns\ds \leq e^{C r_1} \Big(\Big|\frac{\partial
u}{\partial \nu}\Big |_{L^2_{ \mathcal{
F}}(0,T;L^2(\Gamma_0))} + |f|_{L^2_{ \mathcal{
F}}(0,T;H_0^1(G))} + |g|_{L^2_{ \mathcal{
F}}(0,T;H^1(G))}\Big).
\end{array}
\end{equation}}

Unfortunately, following the method in this
paper, one could obtain only an inequality as
follows:
\begin{equation} \label{bobser esti4}
\begin{array}{ll}\ds
 \q |u_T|_{L^2(\Omega,{ \mathcal{F}}_T, P; H_0^1(G))} \\
\ns\ds \leq e^{C r_1} \Big(\Big|\frac{\partial
u}{\partial \nu}\Big |_{L^2_{ \mathcal{
F}}(0,T;L^2(\Gamma_0))} +
|U|_{L^2_\cF(0,T;H^1(G))} + |f|_{L^2_{\mathcal{
F}}(0,T;H_0^1(G))}\\
\ns\ds \qq\qq + |g|_{L^2_{ \mathcal{
F}}(0,T;H^1(G))}\Big).
\end{array}
\end{equation}
It seems to us that  getting rid of  the
undesired term $|U|_{L^2_\cF(0,T;H^1(G))}$ in
the inequality \eqref{bobser esti4} is a very
challenging task.

\item {\bf Construction of the solution $z$ from the observation}

In this paper, we only answer the first and the
second questions in the state observation
problem. The third one is still open. Since the
equation \eqref{system2} is time irreversible,
some efficient approaches (See \cite{Li1} for
example), which work well for time reversible
systems, become invalid. On the other hand, we
may consider the following minimization
problem:

{\it Find a $\bar z_0\in
L^2(\O,\cF_0,P;H_0^1(G))$ such that
$$
\| \frac{\pa \bar z}{\pa\nu} - h
\|_{L^2_\cF(0,T;L^2(\G_0))}=\min_{z_0\in
L^2(\O,\cF_0,P;H_0^1(G))}\| \frac{\pa
z}{\pa\nu} - h \|_{L^2_\cF(0,T;L^2(\G_0))},
$$
where $h\in L^2_\cF(0,T;L^2(\G_0))$ is the
observation and $z$ (\resp$\bar z$) is the
solution to the equation \eqref{system2} with
initial datum $z_0$ (\resp$\bar z_0$).}

\no It seems that one may utilize the method
from optimization theory to study the
construction of $z_0$. Because of the
stochastic nature, this  is an interesting but
difficult problem and the detailed analysis is
beyond the scope of this paper.

\item{\bf Unique continuation property with less restrictive conditions}

In this paper, we  show that, under the
condition \eqref{zx11}, $y=0$ in $Q$, $P$-a.s.
Compared to the classical unique continuation
result for deterministic Schr\"{o}dinger
equations with time independent coefficients
(see \cite{Es1,Lebeau} for example), the
condition \eqref{zx11} is too restrictive. It
would be quite interesting but maybe
challenging to prove whether the result in
\cite{Es1} is true or not for stochastic
Schr\"{o}dinger equations. In fact, as far as
we know, people even do not know whether the
results in \cite{Es1,Lebeau} are true or not
for deterministic Schr\"{o}dinger equations
with time-dependent lower order term
coefficients, which is a particular case of the
equation \eqref{system1}.

\end{itemize}

%%%%%%%%%%%%%%%%%%%%%%%%%%%%%%%%%%%%%%%%%%%%%%%%%%%%%%%%%%%%%%%%%%%%%%%%%%%%%%%%%%%%

\section*{Acknowledgments}

This paper is an improved version of one
chapter of the author's Ph D thesis
(\cite{Luqi2}) accomplished at Sichuan
University under the guidance of Professor Xu
Zhang. The author would like to take this
opportunity to thank him deeply for his help.
The author also highly appreciates the
anonymous referees for their constructive
comments.

\end{document}